\documentclass{siamltex}
\usepackage[utf8]{inputenc}

\usepackage{amssymb}
\usepackage{amsmath}
\usepackage{latexsym}
\usepackage{graphics}
\usepackage{url}
\usepackage{xcolor}
\usepackage{algorithm}
\usepackage{algpseudocode}

\usepackage{tikz,pgfplots}
\pgfplotsset{compat=newest,
every axis/.append style={axis x line=bottom,
                          axis y line=left,
                          scale only axis,
                          y label style={at={(0.0,1.0)},anchor=south west,rotate=-90}
                          },
}

\def\bv{\mathbf{v}}
\def\bx{\mathbf{x}}
\def\by{\mathbf{y}}
\def\bf{\mathbf{f}}
\def\bg{\mathbf{g}}
\def\bi{\mathbf{i}}
\def\b0{\mathbf{0}}

\newcommand{\cC}{\mathcal{C}}

\newcommand{\R}{\mathbb{R}}

\newcommand{\cU}{\mathcal{U}}

\newcommand{\cJ}{\mathcal{J}}

\newcommand{\px}[1]{\partial_{{x} _{#1}}}

\newcommand{\bb}{\mathbf{b}}

\newtheorem{remark}[theorem]{Remark}
\numberwithin{algorithm}{section}

\title{Tensor Decomposition Methods for High-dimensional Hamilton-Jacobi-Bellman Equations}

\author{Sergey~Dolgov\thanks{S. Dolgov is with the Department of Mathematical Sciences,  University of Bath, North Rd, BA2 7AY Bath, United Kingdom ({\tt S.Dolgov@bath.ac.uk}).}
\and Dante~Kalise\thanks{D. Kalise is with the School of Mathematical Sciences, University of Nottingham, University Park Campus, NG7 2QL Nottingham, United Kingdom ({\tt dante.kalise@nottingham.ac.uk}).}
\and Karl~Kunisch\thanks{K. Kunisch is with Institute of Mathematics and Scientific Computing, University of Graz, Heinrichstr. 36, A-8010 Graz, Austria and
	Radon Institute for Computational and Applied Mathematics (RICAM), Altenbergerstra{\ss}e 69, A-4040 Linz, Austria({\tt karl.kunisch@uni-graz.at}).}}

\begin{document}
\maketitle
\begin{abstract}
	A tensor decomposition approach for the solution of high-dimensional, fully nonlinear Hamilton-Jacobi-Bellman equations arising in optimal feedback control of nonlinear dynamics is presented. The method combines a tensor train approximation for the value function together with a Newton-like iterative method for the solution of the resulting nonlinear system. The tensor approximation leads to a polynomial scaling with respect to the dimension, partially circumventing the curse of dimensionality. A convergence analysis for the linear-quadratic optimal control problem is presented. For nonlinear dynamics, the effectiveness of the high-dimensional control synthesis method is assessed in the optimal feedback stabilization of the Allen-Cahn and Fokker-Planck equations with a hundred of variables.
\end{abstract}



\section{Introduction}
Richard Bellman first coined the expression ``curse of dimensionality'' when referring to the overwhelming computational complexity associated to the solution of multi-stage decision processes through dynamic programming, what is nowadays known as Bellman's equation. More than 60 years down the road, the curse of dimensionality has become ubiquitous in different fields such as numerical analysis, compressed sensing and statistical machine learning. However, it is in the computation of optimal feedback policies for the control of dynamical systems where its meaning continues to be most evident. Here, the curse of dimensionality arises since the synthesis of optimal feedback laws by dynamic programming techniques demands the solution of a Hamilton-Jacobi-Bellman (HJB) fully nonlinear Partial Differential Equation (PDE) cast over the state space of the dynamics. This intrinsic relation between the dimensions of the state space of the control system and the domain of the HJB PDE generates computational challenges of formidable complexity even for relatively simple dynamical systems\footnote{As an illustration, consider the simplest double integrator dynamics $\ddot{x}=u$, whose optimal feedback synthesis already requires the solution of a two-dimensional PDE. For a quadrocopter model where the dynamics are described by a 12-dimensional nonlinear dynamical system, the associated HJB PDE has to be solved in $\R^{12}$. We recall that much of the research in computational PDEs is devoted to the solution of problems in physical space, that is $\R^{3+1}$, at most.}.
	Much of the research in control revolves around circumventing this barrier through different trade-offs between dimensionality, performance, and optimality of the control design. Prominent examples of the research landscape shaped by the curse of dimensionality include model order reduction, model predictive control, suboptimal feedback design, reinforcement learning and distributed control. However, the effective computational solution of dynamic programming equations of arbitrarily high dimensions through deterministic methods remains an open quest with fundamental implications in optimal control design. In this paper, we present a computational approach based on tensor decomposition techniques for the solution of high-dimensional HJB PDEs arising in optimal feedback control of systems governed by partial differential equations. We show that for evolution equations arising from the semi-discretization of PDEs, our technique scales at a rate which is at most polynomial (of largest observed order 4) in the dimension. This scaling allows the computation of accurate feedback laws for nonlinear dynamics with over 100 dimensions. We assess our design over a class of challenging problems, including the optimal feedback stabilization of  nonlinear parabolic PDEs such as the Allen-Cahn and Fokker-Planck equations in one and two-spatial dimensions, and in the presence of control constraints.

	\subsection{The Numerical Approximation of HJB PDEs}
	Since the seminal work by Crandall and Lions \cite{CL84}, the approximation of HJB equations through computational PDE methods has been addressed by a range of discretization strategies, most notably finite differences, level-set methods and semi-Lagrangian schemes \cite{FF16}. The aforementioned techniques have proven to overcome the difficulties associated to the fully nonlinear character of the HJB PDEs \cite{NSZ17}. However, they are inherently grid-based schemes suffering from the curse of dimensionality. That is, for a multidimensional ansatz defined from a tensor product of 1d objects, the scaling of the total number of degrees of freedom of the discretization grows exponentially with respect to the dimension of the HJB PDE. In practice, this makes the problem computationally intractable for dimensions larger than 4. This is a fundamental limitation in the context of nonlinear optimal feedback control, where the dimension of the associated HJB PDE is determined by the dimension of the state space of the control system. A partial remedy to this difficulty is the coupling of grid-based discretizations for low-dimensional HJB PDEs with model reduction techniques to lower the dimension of the control system \cite{KVX04,AFV17}. This approach has been successfully applied to dynamics with strong dissipative properties, but its overall performance relies on a good state space sampling and deteriorates for dynamics including transport, delays, or highly nonlinear phenomena \cite{KK14}. While the rigorous design of numerical methods for the solution of very high-dimensional HJB PDEs remains largely an open problem, encouraging results have been obtained over the last years. A non-exhaustive list includes the use of machine learning techniques \cite{SS18,EHJ17,HJE18,HPW19},  approximate dynamic programming in the context of reinforcement learning \cite{B19,RZ19}, causality-free methods and convex optimization \cite{KW16,CDOY}, max-plus algebra methods \cite{E07,AGL08}, polynomial approximation  \cite{KK18,KKK19}, tree structure algorithms \cite{AFS19}, and sparse grids \cite{BGGK13,GK16}. A very recent stream of works \cite{EHJ17,SS18,HPW19,RK18,Hutzenthaler-parabolic-2020} has explored the use of machine learning techniques to approximate high-dimensional nonlinear PDEs. The work \cite{SS18} proposes the so-called Deep Galerkin Method, combining a deep neural network ansatz for the solution together with a PDE residual minimization. In \cite{EHJ17,HJE18,HPW19}, the authors focus on the class of time-dependent HJB equations arising in stochastic control where a pointwise evaluation of the solution can be realized through a representation formula involving the solution of a backward  stochastic differential equation. These latter approaches can be linked to causality-free methods, with deterministic counterparts explored in \cite{KW16,NGK19,CDOY,YD18}. All the aforementioned methods have been assessed in high-dimensional settings and have relevant features which make them appealing for specific instances of the feedback design problem. However, to the best of our knowledge, the methodology proposed in this work is the first PDE-based computational method to go beyond 100 dimensions for the stationary HJB PDE associated to deterministic infinite horizon nonlinear optimal control. While a benchmarking against data-driven methods goes beyond the scope of this paper, we compare our method with the sub-optimal feedback law obtained from solving a linearised control problem, which is a technique widely used in the control literature and which scales well for high-dimensional settings.
	
	\subsection{A tensor calculus framework for nonlinear HJB} In this work we propose a numerical method for the solution of HJB equations based on tensor decomposition techniques \cite{BM02,KB09,hackbusch-2012}, which have proven to be successful in tackling the curse of dimensionality in the context of numerical analysis of PDEs \cite{K15,DV19}. The use of low-rank structures such as the tensor-train (TT) format \cite{osel-tt-2011} to represent high-dimensional objects allows the solution of linear high-dimensional problems by generalizing standard numerical linear algebra techniques to multi-index arrays of coefficients (tensors) and the multivariate functions they approximate. This approach has been recently explored in \cite{HDB14} for solving a class of finite-horizon stochastic control problems where the associated time-dependent HJB PDE can be transformed into a linear equation \cite{TPNAS}.
	A fixed-point iteration algorithm using tensor approximations with a Markov chain discretization was proposed in \cite{Marzouk-hjb-robotics,Marzouk-hjb-CDC}.
	A model-free learning of a TT representation of the value function from Monte Carlo realisations was proposed in \cite{Schneider-hjb-2019}.
	However, all Monte-Carlo-based methods may suffer from slow convergence.
	Here, we extend the tensor calculus framework to approximate the solution of fully nonlinear, first-order, stationary HJB PDEs arising from deterministic infinite horizon control, using a fast deterministic Newton-type policy iteration and a spectral discretization that ensures a fast algebraic (in some examples nearly exponential) convergence, given a high regularity of value functions arising in the considered class of feedback control problems for semi-linear PDEs.
	Furthermore, through the selection of suitable control penalties \cite{Ly98}, our method allows the inclusion of control constraints in the design.
	
	\subsection{Contributions of this work} We develop a computational approach based on tensor decompositions for the solution of Hamilton-Jacobi-Bellman PDEs arising in the computation of optimal actions in feedback form. Our method scales at most polynomially with the state space dimension for a wide class of systems governed by high-dimensional nonlinear evolution equations, including dynamics originating from the Allen-Cahn and Fokker-Planck equations. The mitigation of the curse of dimensionality allows the accurate synthesis of optimal feedback maps for nonlinear dynamics with over 100 dimensions at a moderate computational cost. The method can effectively incorporate control constraints in the design.

The rest of the paper is structured as follows. In Section \ref{shjb}, we formulate the design problem of computing optimal feedback controllers for nonlinear dynamics via the HJB equation, together with an abstract iterative method for its approximation. In Section \ref{stensors} we develop all the building blocks underpinning a tensor decomposition framework for the solution of high-dimensional HJB equations. Finally, in Section \ref{spdes} we apply the proposed methodology to the computation of optimal feedback laws for the Allen-Cahn and Fokker-Planck equations.

	\section{The HJB PDE in optimal feedback control}\label{shjb} We study the following infinite horizon optimal control problem:
	\begin{align}
	\underset{u(\cdot)\in\cU}{\min}\;\cJ(u(\cdot),\bx):=\int\limits_0^\infty \ell(\by(t))+\gamma|u(t)|^2\, dt\,,\label{ihoc}
	\end{align}
	subject to the nonlinear dynamical constraint
	\begin{align}
	\dot \by(t)= \bf(\by(t))+\bg(\by)u(t)\,,\quad \by(0)=\bx,\label{noldyn}
	\end{align}
	where the state $\by(t)=(y_1(t),\ldots,y_d(t))^\top\in\R^d$, the control $u\in\cU\equiv L^{\infty}([0;+\infty[;U)$, with $U$ a compact set of $\R$, the running cost $\ell(\by):\R^d\rightarrow\R_0^+$, and  the control penalization $\gamma>0$. We assume the state cost $\ell(\by)$ and the system dynamics $\bf(\by):\R^d\rightarrow\R^d$ and {$\bg(\by):\R^d\rightarrow\R^{d}$} to be $\cC^1(\R^d)$. Without loss of generality, the origin $\by=\b0$ is an equilibrium of the uncontrolled dynamics and $\ell(\b0)=0$.
	The control problem \eqref{ihoc} corresponds to the design of a globally asymptotically stabilizing control signal $u(t)$, which can be solved by dynamic programming techniques. Defining the value function
	\begin{align}
	V(\bx):=\underset{u(\cdot)\in\cU}{\inf} \cJ(u(\cdot),\bx)\,,
	\end{align}
	we characterize the solution of the infinite horizon control problem as the unique viscosity solution of the HJB PDE
	\begin{align}\label{hjb}
	\underset{u\in U}{\min}\{ (\bf(\bx)+\bg(\bx) u)^\top DV(\bx)+ \ell(\bx)+\gamma|u|^2\}=0\,,
	\end{align}
	where $D V(\bx)=(\px{1}V,\ldots,\px{d}V)^\top$.  In the unconstrained case $U\equiv \R$, the minimizer $u^*$  is expressed in a feedback form
	\begin{equation}\label{optc}
	u^*(\bx)=-\frac{1}{2\gamma} \bg(\bx)^\top DV(\bx)\,,
	\end{equation}
	which after inserting in \eqref{hjb} leads to the HJB equation
	\begin{equation}\label{hjb2}
	DV(\bx)^\top \bf(\bx)-\frac{1}{4\gamma}DV(\bx)^\top \bg(\bx)\bg(\bx)^\top DV(\bx)+\ell(\bx)=0\,.
	\end{equation}
	This derivation can be extended to controls in $\R^m$ with $m>1$. Moreover, it can be
	modified to enforce box constraints in the control action \cite{Ly98}, by replacing the control penalty $\gamma|u|^2$ by
	\begin{equation}
	W(u)=2\gamma\int\limits_0^u \mathcal{P}^{-1}(\mu)d\mu\,,
	\end{equation}
	where $\mathcal{P}:\R\to\R$ is an odd, bounded, integrable, bijective $\mathcal{C}^1$ function.
	The optimal feedback is given by
	\begin{equation}\label{optc2}
	u^*(\bx)=-\mathcal{P}\left(\frac{1}{2\gamma} \bg(\bx)^\top DV(\bx)\right)\,,
	\end{equation}
	where we can impose lower and upper bound constraints by choosing penalties of the type $\mathcal{P}(x)=u_{\max}\cdot\tanh(x/u_{\max})$. Note that such a choice for the control penalty generates a signal which effectively remains inside $[-u_{\max},u_{\max}]$, however, it differs from the optimal signal that would be obtained by casting a hard control constraint and replacing \eqref{optc2} by the projection
	\begin{equation}
		u^*(\bx)=\mathbb{P}_{[-u_{\max},u_{\max}]}\left(-\frac{1}{2\gamma} \bg(\bx)^\top DV(\bx)\right)\,.
	\end{equation}
	
	\subsection{An iterative approach for solving nonlinear HJB PDEs}
	The construction of a numerical scheme for \eqref{hjb} begins by dealing with the quadratic nonlinearity in the gradient. We apply the  Continuous Policy Iteration developed in \cite{BST7}, a variant of the well-known policy iteration algorithm in dynamic programming \cite{B55,PB79,AFK15}. Conceptually speaking, given an initial guess $u_0(\bx)$ for the optimal feedback control, we insert it into \eqref{hjb} which then becomes a linear PDE for $V(\bx)$, whose solution dictates the update of the feedback control via \eqref{optc}. Algorithm \ref{alg:sga1} summarizes the main steps of the procedure. The algorithm is equivalent to the application of a Newton-type method for $V(\bx)$ directly over \eqref{hjb2}. To guarantee the convergence of the policy iteration when solving \eqref{hjb} over a subdomain $\Omega\subset\R^d$, we require an initial feedback map $u_0(\bx)$ such that $\cJ(u_0(\bx),\bx)<\infty$ over $\Omega$.
	\begin{algorithm}[H]	
		\begin{algorithmic}[1]
			\Require Admissible feedback $u_0(\bx)$, stopping tolerance $\delta>0$
			\While{$error>\delta$}
			\State \textbf{Solve the linearized HJB:} $(\bf(\bx)+\bg(\bx)u_s)^TDV_s(\bx)+ \ell(\bx)+\gamma|u_s|^2=0\,\,.$
			\State \textbf{Feedback update:}  $u_{s+1}(\bx)=-\mathcal{P}\left(\frac{1}{2\gamma}\bg^\top DV_s(\bx)\right)$.
			\State \textbf{Set } $error=\|V_{s}-V_{s-1}\|$, $s:=s+1$.
			\EndWhile \\
			\Return $(V_s,u_s)\approx(V^*,u^*)$
		\end{algorithmic}
		\caption{Continuous Policy Iteration Algorithm}\label{alg:sga1}
	\end{algorithm}	
	
	\section{A tensor decomposition framework for high-dimensional HJB equations}\label{stensors}
	We employ the Galerkin spectral element approximation similarly to \cite{KK18} except that we construct the basis functions from the Legendre polynomials of bounded maximal \emph{individual} degree $n-1$,
	\begin{equation}
	\Phi_{\bi}(\bx):=\phi_{i_1}(x_1)\cdots\phi_{i_d}(x_d), \quad i_k=0,\ldots,n-1,
	\label{eq:Vnspace}
	\end{equation}
	where $\phi_{i_k}(x_k)$ are the univariate Legendre polynomials of degree $i_k\leq n-1$, and the multi-indices $\bi = (i_1,\ldots,i_d)$, $\mathbf{j} = (j_1,\ldots,j_d)$.
	In the $s$th iteration of Alg.~\ref{alg:sga1}, we seek the value function in the form
	\begin{equation}
	V_s(x_1,\ldots,x_d) \approx \sum_{j_1,\ldots,j_d=0}^{n-1} \mathbf{v}(j_1,\ldots,j_d) \Phi_{j_1,\ldots,j_d}(\bx),
	\label{eq:ansatz}
	\end{equation}
	by making the Galerkin residual of the linearized HJB orthogonal to $\{\Phi_{\bi}\}$.
	This requires solving a system of $n^d$ Galerkin equations in $n^d$ unknowns $\mathbf{v}(\mathbf{j})$,
	\begin{equation}
	\sum_{\mathbf{j}} \underbrace{\left\langle \Phi_{\bi}, (\bf+\bg u_s)^\top D \Phi_{\mathbf{j}}\right\rangle}_{A(\mathbf{i},\mathbf{j})} \bv(\mathbf{j})  = \underbrace{-\left\langle \Phi_{\bi}, ~\ell + \gamma u_s^2 \right\rangle}_{\bb(\mathbf{i})},
	\label{eq:discnd}
	\end{equation}
	where $\langle \cdot, \cdot \rangle$ is an inner product in $L^2\left([-a,a]^d\right)$ with an appropriately chosen \emph{domain size} $a>0$. Given the tensor product structure of \eqref{eq:ansatz}, its accuracy can be analyzed with univariate polynomial approximation theory~\cite{trefethen-spectral-2000}, with an exponential error decay rate $\mathcal{O}(n^{-p})$ for $V(\bx)\in\cC^p(\Omega)$ (e.g. $\cC^{\infty}(\mathbb{R}^d)$, see Remark~\ref{rem:V-reg} below).

	\subsection{Compressed Tensor Train representation}
	The coefficients in \eqref{eq:ansatz} are enumerated by $d$ independent indices, so $\mathbf{v}$ can be treated as a $d$-dimensional \emph{tensor}.
	Throughout the paper, we approximate such tensors by the so-called Tensor Train (TT) decomposition \cite{osel-tt-2011},
	\begin{equation}
	\mathbf{\tilde v}(\mathbf{i}) := \sum_{\alpha_0,\ldots,\alpha_d=1}^{r_0,\ldots,r_d} \bv^{(1)}_{\alpha_0,\alpha_1}(i_1) \bv^{(2)}_{\alpha_1,\alpha_2}(i_2)  \cdots \bv^{(d)}_{\alpha_{d-1},\alpha_d}(i_d).
	\label{eq:tt}
	\end{equation}
        The smaller (3-dimensional) tensors $\bv^{(k)}$ on the right hand side are called \emph{TT blocks},
        and the new summation ranges $r_0,\ldots,r_d$ are called \emph{TT ranks}.
        Notice that when $r_0=r_1=\cdots=r_d=1$ the complete \emph{separation of variables} is attained, as the tensor becomes a product of univariate factors.
        In general we can consider $r_1,\ldots,r_{d-1}$ to be larger than $1$.
        In case of two variables, the TT decomposition becomes just a dyadic factorisation of a low-rank matrix.
        The optimal low-rank matrix approximation can be computed using the singular value decomposition (SVD),
        and a general TT decomposition can be computed for any tensor by $d-1$ SVDs with quasi-optimal TT ranks for the given accuracy $\varepsilon$ \cite{osel-tt-2011}.
        For convenience we can introduce the maximal TT rank $r := \max_{k=0,\ldots,d} r_k$.
        Counting the number of unknowns in the TT blocks in \eqref{eq:tt}, one can conclude that the TT decomposition needs at most $dnr^2$ unknowns.
        For numerical efficiency, we assume that $r$ can be taken much smaller than the original cardinality $n^d$ for a desired approximation accuracy.
        This assumption is fulfilled for many relevant cases, as we show next.

        For theoretical analysis, it is convenient to combine \eqref{eq:tt} with \eqref{eq:ansatz} and to consider the \emph{functional} TT format \cite{osel-constr-2013},
        $$
        V(\bx) \approx \widetilde V(\bx) := \sum_{\alpha_0,\ldots,\alpha_d=1}^{r_0,\ldots,r_d} v^{(1)}_{\alpha_0,\alpha_1}(x_1) \cdots v^{(d)}_{\alpha_{d-1},\alpha_d}(x_d),
        $$
        to work directly over the TT ranks of a function.
        Sharp TT rank bounds are usually hard to derive though: the SVD might reveal an optimal decomposition that is difficult to express analytically.
        It was proven for certain (such as rational) classes of functions \cite{tee-tensor-2003,uschmajew-approx-rate-2013,grasedyck-kron-2004} and extensively tested numerically in wider scenarios
        that smooth (e.g. analytic) functions exhibit a logarithmic growth of TT ranks, $r \sim \log^p (1/\varepsilon)$, to achieve an error $\varepsilon$.
        As a rationale for using the TT format for the HJB equations,
        here we show that quadratic value functions admit TT approximations with a similar poly-logarithmic convergence rate.
	
	\begin{theorem}\label{thm:rank}
	Assume that $\ell(\by) = \frac{1}{2}\by^T Q \by$, where $Q$ is a symmetric positive definite matrix  with the Cholesky decomposition $Q=D^\top D$, and that the linear system $\dot \by(t) = A\by(t) + Bu$ is stabilisable. There exists a solution $\Pi\in\mathbb{R}^{d \times d}$ of the Riccati equation
	$$
	A^\top \Pi + \Pi A - \frac{1}{\gamma} \Pi BB^\top \Pi + Q = 0,
	$$
	such that the eigenvalues of $A_{\pi} = AD^{-1} - \frac{1}{\gamma} BB^\top\Pi D^{-1}$ satisfy $\lambda(A_{\pi}) \in [\lambda_{\min}, \lambda_{\max}] \oplus \mathrm{i}[-\mu,\mu],$ $\lambda_{\max}<0$.  Assume that the ranks of the off-diagonal blocks of $A_{D}=AD^{-1}$ are bounded by a constant, $\mathrm{rank}~A_{D}(k+1:d,~1:k) \le M$ for all $k=1,\ldots,d-1$,
	and that $\mathrm{rank}(B) \le r_b$.
	Then for any $\varepsilon \in (0,1)$ the value function $V(\bx) = \bx^\top \Pi \bx$ admits a TT approximation $\widetilde V(\bx)$ with the TT ranks
	$$
	r_k \le \min\left((M+r_b)\left(\log\frac{1}{\varepsilon} + C\right)^{7/2},~ \min(k,d-k)\right)+2,
	$$
	and the error $\max_{\bx \in [-a,a]^d} |V(\bx) - \widetilde V(\bx)| \le \varepsilon$
	for some offset $$C = C_0 + \frac{\mu}{|\lambda_{\max}|} + 2\log \left[\frac{\sqrt{\lambda_{\min}^2+\mu^2}}{|\lambda_{\max}|}\frac{a\|A_{\pi}\|\|D^{-\top}\Pi B\|}{\gamma}\right]>0\,,$$ where $C_0$ is independent of $d,\varepsilon,M,r_b,\gamma,\mu,\lambda_{\min},\lambda_{\max}$.
	If the second bound $r_k = \min(k,d-k)+2$ is attained for all $k$, the TT decomposition $\widetilde V$ is exact.
	\end{theorem}
        \begin{proof}
	For the upper bound we employ \cite[Thm 4.2]{dk-qtt-tucker-2013}: for the second order polynomial $V(\bx) = \bx^\top \Pi \bx$
	with a symmetric matrix $\Pi$
	there exists an exact TT decomposition
	with the TT ranks governed by the ranks of the off-diagonal blocks of $\Pi$,
	\begin{equation}
	r_k \le \mathrm{rank}\left(\Pi(k+1:d,~1:k)\right)+2.
	\label{eq:rank-xPx}
	\end{equation}
	This gives an obvious bound $\min(k,d-k)+2$.
	
	However, there might exist an approximate TT decomposition of lower TT ranks.
	First, we notice that the Riccati equation can be rewritten as a Lyapunov equation.
	In fact, using a stabilised matrix \cite{bkp-ns-control-2019} $A - \frac{1}{\gamma} BB^\top\Pi$, and also the Cholesky factor of $Q$, we obtain
	\begin{equation}
	A_{\pi}^\top \Pi + \Pi A_{\pi} =  -\frac{1}{\gamma} D^{-\top} \Pi B B^\top \Pi D^{-1} - I,
	\label{eq:Alyap}
	\end{equation}
	where $I$ is a $d \times d$ identity matrix.
	The left hand side is constructed from the stable matrix $A_{\pi}$.
	Using the Kronecker product, we can write the Lyapunov equation as a large linear system,
	$$
	\underbrace{\left(A_{\pi} \otimes I + I \otimes A_{\pi}\right)}_{\mathcal{A}} \mathrm{vec}(\Pi) = -\frac{1}{\gamma} \sum_{i=1}^{r_b} (D^{-1} \Pi B_i) \otimes (D^{-1} \Pi B_i) - \mathrm{vec}(I),
	$$
	where $\mathrm{vec}(\cdot)$ stacks all columns of a matrix into a vector.
	Now we can use \cite[Thm. 9]{grasedyck-kron-2004}: there exists an approximate inverse $\mathcal{\widetilde A}^{-1}$ of the Kronecker product form
	$$
	\mathcal{\widetilde A}^{-1} = \sum_{j=-R}^{R} \frac{2 w_j}{\lambda_{\max}} \exp\left(-\frac{2t_j}{\lambda_{\max}} A_{\pi}\right) \otimes \exp\left(-\frac{2t_j}{\lambda_{\max}} A_{\pi}\right)
	$$
	with the approximation error
	\begin{equation}
	\|\mathcal{A}^{-1} - \mathcal{\widetilde A}^{-1}\| \le \tilde C \|\mathcal{A}\|\frac{\sqrt{\lambda_{\min}^2+\mu^2}}{|\lambda_{\max}|}  \exp\left(\frac{2}{\pi}\frac{\mu}{|\lambda_{max}|} - \pi \sqrt{2R}\right),
	\label{eq:errR}
	\end{equation}
	where $\tilde C>0$ is a constant independent of other parameters in~\eqref{eq:errR}.
	Multiplying the approximation $\mathcal{\widetilde A}^{-1}$ with the low-rank first term in the right hand side of \eqref{eq:Alyap},
	we obtain an incomplete solution of the form
	$\mathrm{vec}(\widehat\Pi) = \sum_{i=1}^{(2R+1)r_b} p_i \otimes q_i$,
	and hence the rank of $\widehat\Pi$ is bounded by $(2R+1)r_b$.
	The negative identity matrix in \eqref{eq:Alyap} yields the second term
	\begin{equation}
	\mathrm{vec}(\check\Pi) = \mathcal{\widetilde A}^{-1} \mathrm{vec}(-I), \quad \check\Pi= \sum_{j=-R}^{R} -\frac{2 w_j}{\lambda_{\max}} \exp\left(-\frac{2t_j}{\lambda_{\max}} (A_{\pi}+A_{\pi}^\top)\right)
	\label{eq:Pi2}
	\end{equation}
	in the ultimate approximate solution $\widetilde\Pi = \widehat\Pi + \check\Pi$.
	
	Since the first term $\widehat\Pi$ is a low-rank matrix, all its off-diagonal blocks \eqref{eq:rank-xPx} have low ranks of at most $(2R+1)r_b$ too.
	However, $\check\Pi$ is a full-rank matrix and we need to investigate its off-diagonal, also called \emph{quasi-separable} \cite{Eidelman-qsep-1999}, ranks directly.
	First, we recall \cite[Lemma 16]{grasedyck-kron-2004} that each matrix exponential in \eqref{eq:Pi2} can be approximated by a sum of $2k_e+1$ resolvents with an error
	\begin{equation}
	\begin{split}
	& \left\|\exp\left(-\frac{2t_j}{\lambda_{\max}} (A_{\pi}+A_{\pi}^\top)\right) - \sum_{\ell=-k_e}^{k_e}\kappa_{\ell} \left(z_{\ell}I + \frac{2t_j}{\lambda_{\max}} (A_{\pi}+A_{\pi}^\top)\right)^{-1} \right\| \\
	& \le \bar C \exp\left(4\left(\frac{4t_j \mu}{|\lambda_{max}|}+1\right)^2 - \left(\frac{4 t_j \mu}{|\lambda_{max}|}+1\right)^{2/3} k_e^{2/3}\right).
	\end{split}
	\label{eq:errk}
	\end{equation}
	Since the quasi-separable rank of $A_{z_{\ell}}:=z_{\ell}I + \frac{2t_j}{\lambda_{\max}} (A_{\pi}+A_{\pi}^\top)$ coincides with that of $A_{\pi}$, which is $M+r_b$,
	and on the other hand it coincides with the quasi-separable rank of the inverse matrix $A_{z_{\ell}}^{-1}$ \cite{Eidelman-qsep-1999,Hackbusch-H-1999},
	we can conclude that the approximate quasi-separable rank of $\exp\left(-\frac{2t_j}{\lambda_{\max}} (A_{\pi}+A_{\pi}^\top)\right)$
	is bounded by $(2k_e+1)(M+r_b)$,
	and the quasi-separable rank of $\check\Pi$ \eqref{eq:Pi2} is bounded by $(2R+1)(2k_e+1)(M+r_b)$.
	
	The approximate value function is constructed as $\widetilde V(\bx) = \bx^\top \widetilde\Pi \bx$,
	and by \eqref{eq:rank-xPx} we can estimate its TT rank as
	\begin{equation}
	r_k(\widetilde V) \le (2R+1)(r_b + (2k_e+1)(M+r_b)) + 2.
	\label{eq:rank-raw}
	\end{equation}
	
	For the error estimate, we have
	$$
	\varepsilon = \max_{\bx \in [-a,a]^d} |V(\bx) - \widetilde V(\bx)| \le a^2 \|\Pi - \widetilde\Pi\| \le a^2 \|\mathcal{A}^{-1} - \mathcal{\widetilde A}^{-1}\| \left(\frac{\|D^{-\top} \Pi B\|^2}{\gamma}+1\right).
	$$
	From \eqref{eq:errR} we obtain
	\begin{equation*}
	\begin{split}
	R & \le \frac{1}{2\pi^2}\left(\frac{2}{\pi}\frac{\mu}{|\lambda_{\max}|} + \log \|A\| + \log\frac{\sqrt{\lambda_{\min}^2+\mu^2}}{|\lambda_{\max}|} + \hat C + \left|\log\|\mathcal{A}^{-1} - \mathcal{\widetilde A}^{-1}\|\right|\right)^2, \\
	& \le \left(\log\frac{1}{\varepsilon} + \frac{\mu}{|\lambda_{\max}|} + \log\left(\frac{a^2 \|D^{-\top} \Pi B\|^2 \|A_{\pi}\| \sqrt{\lambda_{\min}^2+\mu^2}}{\gamma |\lambda_{\max}|} \right) + \hat C\right)^2
	\end{split}
	\end{equation*}
	for some constant $\hat C>0$,
	while \eqref{eq:errk} together with \cite[Lemma 5]{grasedyck-kron-2004} gives
	$$
	k_e \le \left(\left|\log \|\mathcal{A}^{-1} - \mathcal{\widetilde A}^{-1}\|\right| + \log R + \check C\right)^{3/2}
	$$
	for $\check C>0$ being some other constant.
	Plugging these bounds into \eqref{eq:rank-raw}, we obtain the first estimate of the TT rank.

	\hfill\end{proof}
	
	
	\begin{remark}
		In many cases one can take $Q$, and hence $D$, to be diagonal matrices, for example, if the 2-norm of the state vector (corresponding to the $L^2$-norm of the state function) is used in the cost functional.
		In this case the ranks of the off-diagonal blocks of $A_D$ coincide with those of $A$.
	\end{remark}
	
	\begin{remark}
		Although Thm.~\ref{thm:rank} is formulated only for linear systems,
		the remarkable proportionality between the TT ranks of the value function and the off-diagonal ranks of the linearised system matrix
		seems to hold for some semi-linear systems as well.
		In particular, we observe from Fig.~\ref{fig:nw-d} that the TT ranks of the value function for discretised one-dimensional PDEs grow very mildly with the number of variables,
		which can be also attributed to the growth of the ratios $\lambda_{min}/\lambda_{max}$, $\mu/\lambda_{\max}$.
		However, when the system is produced from a two-dimensional PDE, the TT ranks grow proportionally to the number of degrees of freedom introduced in each spatial direction (see Fig.~\ref{fig:nw2-d}).
		Thus, the ranks of the actuator matrix and of the off-diagonal blocks of the Jacobian matrix
		can give a useful hint whether the TT approach might be efficient for the HJB equation of a particular dynamical system of interest.
	\end{remark}
	
	\begin{remark}\label{rem:V-reg}
		Alternatively, fast convergence of the TT approximation can be related to the smoothness of the original function~\cite{tee-tensor-2003,uschmajew-approx-rate-2013}.
		For example, it was verified in ~\cite{BreKP18} that the cost functional for the bilinear optimal control problem for the Fokker-Planck equation we present in our numerical results belongs to~$\cC^{\infty}$ in a neighborhood of the origin. In \cite{kw2020} sufficient conditions for $\cC^{1}$-regularity of the value function
on bounded sets are given. 
	\end{remark}
	
	\subsection{TT decomposition of the system functions and HJB equation}
        To take advantage of the TT decomposition of the value function, it is necessary to find also compatible representations of the stiffness matrix and right hand side of the Galerkin HJB equations~\eqref{eq:discnd}.
        For example,
	matrices of size $n^d \times n^d$, such as that in the left hand side of~\eqref{eq:discnd}, can be represented in a slightly different \emph{matrix} TT format, where we separate pairs of row and column indices,
        \begin{equation}
        A(\bi,\mathbf{j}) = \sum_{\beta_0,\ldots,\beta_d=1}^{R_0,\ldots,R_d} A^{(1)}_{\beta_0,\beta_1}(i_1,j_1) \cdots A^{(d)}_{\beta_{d-1},\beta_d}(i_d,j_d).
        \label{eq:ttm}
        \end{equation}
        For complexity estimates we define also the upper bound $R:=\max_{k=0,\ldots,d} R_k$.
        In particular, we need to construct linear parts of the stiffness matrix,
        \begin{equation}\label{eq:Alin}
        A_{f_p}(\bi,\mathbf{j}) := \langle \Phi_{\bi}, f_p \px{p} \Phi_{\mathbf{j}} \rangle = \int_{[-a,a]^d} \Phi_{\bi}(\bx) f_p(\bx) \px{p}\Phi_{\mathbf{j}}(\bx) d\bx,
        \end{equation}
        where $f_p(\bx)$ is the $p$-th component of the drift $\bf(\bx) = (f_1(\bx),\ldots,f_d(\bx))$, $p=1,\ldots,d$.
        The integral in \eqref{eq:Alin} can be approximated by a tensorised Gauss-Legendre quadrature, e.g.
        $$
        A_{f_p}(\bi,\mathbf{j}) \approx \sum_{k_1,\ldots,k_d=1}^{m} w_{k_1}\cdots w_{k_d} f_p(x_{k_1},\ldots,x_{k_d}) \Phi_{\bi}(x_{k_1},\ldots,x_{k_d}) \px{p}\Phi_{\mathbf{j}}(x_{k_1},\ldots,x_{k_d}),
        $$
        where $x_k,w_k$, $k=1,\ldots,m$, are quadrature nodes on $(-a,a)$ and weights, respectively.
        For efficiency we always take $m=\mathcal{O}(n)$, e.g. $m=2n$.
        Suppose that a TT decomposition of $f_p(x_{k_1},\ldots,x_{k_d})$ is given,
        \begin{equation}
        \bf_p(k_1,\ldots,k_d) := f_p(x_{k_1},\ldots,x_{k_d}) \approx \sum_{\beta_0,\ldots,\beta_d=1}^{R_0,\ldots,R_d} \bf^{(1)}_{p,\beta_0,\beta_1}(k_1) \cdots \bf^{(d)}_{p,\beta_{d-1},\beta_d}(k_d).
        \label{eq:fsampled}
        \end{equation}
        Plugging \eqref{eq:fsampled} into \eqref{eq:Alin} and distributing the summations, we can compute directly the TT blocks of a TT decomposition of the stiffness matrix part,
        $$
        A_{f_p}(\bi,\mathbf{j}) = \sum_{\beta_0,\ldots,\beta_d=1}^{R_0,\ldots,R_d}
        A^{(1)}_{p,\beta_0,\beta_1}(i_1,j_1) \cdots A^{(d)}_{p,\beta_{d-1},\beta_d}(i_d,j_d),
        $$
        where
        \begin{equation}
        A^{(p)}_{p,\beta_{p-1},\beta_p}(i_p,j_p) = \left(\sum_{k_p=1}^{m} w_{k_p} \phi_{i_p}(x_{k_p}) \bf^{(p)}_{p,\beta_{p-1},\beta_p}(k_p) \frac{d\phi_{j_p}(x_{k_p})}{dx} \right),
        \label{eq:ttm_fpp}
        \end{equation}
        and
        \begin{equation}
        A^{(q)}_{p,\beta_{q-1},\beta_q}(i_q,j_q) = \left(\sum_{k_q=1}^{m} w_{k_q} \phi_{i_q}(x_{k_q}) \bf^{(q)}_{p,\beta_{q-1},\beta_q}(k_q) \phi_{j_q}(x_{k_q}) \right)
        \label{eq:ttm_fpq}
        \end{equation}
        for $q \neq p$.
        Summing all different components $A_{f_p}$ over $p=1,\ldots,d$, one obtains the complete linear part $\langle \Phi_{\bi}, \bf^\top D \Phi_{\mathbf{j}} \rangle_{L^2(\Omega)}$.
        This summation can be performed in the TT format directly, followed by a TT rank truncation \cite{osel-tt-2011}, or using the iterative Alternating Linear Scheme approximation (see Sec.~\ref{sec:iteration} and \cite{holtz-ALS-DMRG-2012}).
        In our numerical calculations we use the latter approach which requires $\mathcal{O}(d^2 n^2 R^2 r^2)$ floating point operations.
        Similarly we can compute the right hand side entries $\mathbf{b}(\bi) = \langle -\Phi_{\bi}(\bx), \ell(\bx) + \gamma u_s(\bx)^2 \rangle_{L^2(\Omega)}$,
        as well as nonlinear parts of the stiffness matrix, provided that tensors of nodal values $\bg(x_{k_1},\ldots,x_{k_d})$ and $u_s(x_{k_1},\ldots,x_{k_d})$ are approximated by TT decompositions, similarly to~\eqref{eq:fsampled}.

	The system functions are approximated in the TT format \eqref{eq:fsampled}
	using another iterative procedure, the so-called \emph{TT-Cross} algorithm~\cite{ot-ttcross-2010}.
	Here we only recall the main idea of the TT-Cross method to illustrate the key operations needed.
	Any exact TT decomposition, e.g. in the form \eqref{eq:fsampled}, can be recovered from samples of the original tensor by an \emph{interpolation formula} \cite{sav-qott-2014}
	\begin{equation}
	\bf_p(\bi) = \sum_{\substack{\beta_k,\beta_k'=1 \\ k=1,\ldots,d-1}}^{R_k}\bf_p(i_1,\mathcal{I}^{>1}_{\beta_1'}) \left(\bf_p(\mathcal{I}^{\le 1},\mathcal{I}^{>1})\right)^{-1}_{\beta_1',\beta_1} \bf_p(\mathcal{I}^{\le 1}_{\beta_1},i_2,\mathcal{I}^{>2}_{\beta_2'}) \cdots \bf_p(\mathcal{I}^{\le d-1}_{\beta_{d-1}}, i_d),
	\label{eq:ttinterp}
	\end{equation}
	where $\mathcal{I}^{\le k} = \{i_1^{\beta_k},\ldots,i_k^{\beta_k}\}_{\beta_k=1}^{R_k}$ and $\mathcal{I}^{>k} = \{i_{k+1}^{\beta_k},\ldots,i_d^{\beta_k}\}_{\beta_k=1}^{R_k}$ are \emph{left}, respectively, \emph{right} index sets chosen such that the intersection matrices $\bf_p(\mathcal{I}^{\le k}, \mathcal{I}^{>k})$ are invertible.
	For uniformity of notation, we let $\mathcal{I}^{\le 0} = \mathcal{I}^{>d} = \emptyset$.
	Note that the inverse intersection matrices can be multiplied with the adjacent three-dimensional factors to obtain the TT blocks of~\eqref{eq:fsampled}, e.g.
        $$
        \bf^{(k)}_{p,\beta_{k-1},\beta_k}(i_k) = \sum_{\beta_k'=1}^{R_k} \bf_p(\mathcal{I}^{\le k-1}_{\beta_{k-1}},i_k,\mathcal{I}^{>k}_{\beta_k'}) \left(\bf_p(\mathcal{I}^{\le k},\mathcal{I}^{>k})\right)^{-1}_{\beta_k',\beta_k}\,.
        $$
	For numerical stability, the inversion is computed via the QR decomposition \cite{ot-ttcross-2010} or the incremental LU decomposition \cite{sav-qott-2014}.
	
	In practice, one seeks an approximate decomposition of the form \eqref{eq:ttinterp}.
	In this case it becomes important to find indices that not only give invertible intersection matrices, but deliver a small approximation error.
	The TT-Cross algorithm optimises the index positions iteratively.
	In the first step, assume that the right sets $\mathcal{I}^{>k}$ are given (for example at random).
	One can compute the first factor $F^{\{1\}}(i_1,\beta_1) = \bf_p(i_1,\mathcal{I}^{>1}_{\beta_1})$, which can be seen as an $m \times R_1$ matrix.
	The smallest approximation error among sampling rank-$R_1$ approximations
	is given by such $i_1 \in \mathcal{I}^{\le 1} \subset \{1,\ldots,m\}$ that select the submatrix of maximum volume, i.e. $|\det F^{\{1\}}(\mathcal{I}^{\le 1},:)| = \max_{\#\mathcal{I}=R_1} |\det F^{\{1\}}(\mathcal{I},:)|$.
	This set can be found by the \emph{maxvol} algorithm \cite{gostz-maxvol-2010} in $\mathcal{O}(mR_1^2)$ operations, similarly to the LU decomposition with pivoting.
	
	Assume now that we have the left index set $\mathcal{I}^{\le k-1}$ and the right set $\mathcal{I}^{>k}$.
	We can compute $R_{k-1} m R_k$ elements of the tensor and arrange them as a $R_{k-1}m \times R_k$ matrix with elements
	$F^{\{k\}}(\beta_{k-1} i_k,~\beta_k) = \bf_p(\mathcal{I}^{\le k-1}_{\beta_{k-1}},i_k,\mathcal{I}^{>k}_{\beta_k})$.
	Now we can apply the maxvol algorithm to $F^{\{k\}}$ to derive the next index set $\mathcal{I}^{\le k}$ as a subset of the union of $\mathcal{I}^{\le k-1}$ and $i_k$.
	This recursive procedure continues until the last TT block $\bf_p(\mathcal{I}^{\le d-1}, i_d)$ is computed.
	Moreover, we can reverse it in a similar fashion and carry out several TT-Cross iterations, as shown in Algorithm \ref{alg:ttcross}.
	This allows to optimise all index sets and, consequently, the approximations \eqref{eq:ttinterp}, \eqref{eq:fsampled} even if the initial guess was inaccurate.
	
        \begin{algorithm}[h!]
                \caption{TT-Cross~\cite{ot-ttcross-2010} iteration for the TT approximation~\eqref{eq:fsampled}}
                \label{alg:ttcross}
                \begin{algorithmic}[1]
                        \Require Initial sets $\mathcal{I}^{>k} \in \mathbb{N}^{R_{k} \times d-k}$, $k=1,\ldots,d-1.$
                        \For {$k=1,2,\ldots,d$}
                        \State Evaluate $R_{k-1}mR_k$ elements $F^{\{k\}}(\beta_{k-1} i_k;~\beta_k) := \bf_p\left(\mathcal{I}^{\le k-1}_{\beta_{k-1}}, i_k, \mathcal{I}^{>k}_{\beta_k}\right).$
                        \State Apply \emph{maxvol} algorithm to $F^{\{k\}}$ to obtain $\mathcal{I}^{\le k} \subset \mathcal{I}^{\le k-1} \cup \{i_k\}$.
                        \EndFor
                        \For {$k=d,d-1,\ldots,2$}
                        \State Evaluate $R_{k-1}mR_k$ elements $F^{\{k\}}(\beta_{k-1};~i_k\beta_k) := \bf_p\left(\mathcal{I}^{\le k-1}_{\beta_{k-1}}, i_k, \mathcal{I}^{>k}_{\beta_k}\right).$
                        \State Apply \emph{maxvol} algorithm to $(F^{\{k\}})^\top$ to obtain $\mathcal{I}^{>k-1} \subset \{i_k\} \cup \mathcal{I}^{>k}$.
                        \EndFor
                        \State Assemble TT blocks $\bf^{(k)}_{p,\beta_{k-1},\beta_k}(i_k) = \sum_{\beta_k'} \bf_p(\mathcal{I}^{\le k-1}_{\beta_{k-1}},i_k,\mathcal{I}^{>k}_{\beta_k'}) \left(\bf_p(\mathcal{I}^{\le k},\mathcal{I}^{>k})\right)^{-1}_{\beta_k',\beta_k}.$
                \end{algorithmic}
        \end{algorithm}

	\begin{remark}
		The result of TT-Cross might still depend on the heuristically chosen initial indices.
		Therefore, we distinguish the stopping threshold (called $\delta$ from now on) and the actual approximation error $\varepsilon$ in the rest of the paper.
	\end{remark}

	Having computed TT approximations to all components $f_p$, we construct the matrix TT blocks \eqref{eq:ttm_fpp}--\eqref{eq:ttm_fpq}.
	Similarly, we apply the TT-Cross algorithm to construct all components of $\bg u_s$, assemble the corresponding parts of the stiffness matrix \eqref{eq:discnd} in the TT format, and sum them together.
	
        We can also precompute a TT matrix of the form \eqref{eq:ttm} which maps the value function coefficients into corresponding control values $\mathbf{u}(\mathbf{j}):=u_{s+1}(x_{j_1},\ldots,x_{j_d})$ on the quadrature grid.
        In the unconstrained control case $u_{s+1} = -\frac{1}{2\gamma} \bg^\top D V_s$, and hence we assemble
	\begin{equation}
	\hat B(\mathbf{j}, \bi) = \sum_{p=1}^{d}\sum_{\beta_0,\ldots,\beta_d} -\frac{1}{2\gamma} \left(\bg^{(1)}_{p,\beta_0,\beta_1}(j_1) \phi^{\delta(p,1)}_{i_1}(x_{j_1}) \right) \cdots \left(\bg^{(d)}_{p,\beta_{d-1},\beta_d}(j_d) \phi^{\delta(p,d)}_{i_d}(x_{j_d}) \right),
	\label{eq:Bcontrol}
	\end{equation}
	using the TT blocks $\bg^{(k)}_{p,\beta_{k-1},\beta_k}(j_k)$ of TT approximations of $g_p(\bx)$, where
	$$
	\phi_i^{\delta(p,q)} = \left\{\begin{array}{ll}d\phi_i/dx, & p=q, \\ \phi_i, & \mbox{otherwise.} \end{array}\right.
	$$
	Now the control signal $\mathbf{u} \approx \hat B \bv$ can be constructed simply as a sum of products of a TT matrix \eqref{eq:ttm} and a TT tensor \eqref{eq:tt}, again with $\mathcal{O}(d^2 n^2 R^2 r^2)$ complexity \cite{osel-tt-2011}.
	In the constrained control case, the first step is the same, followed by approximating the pointwise constraint function
	\begin{equation}
	\mathbf{u}(j_1,\ldots,j_d) = \tilde u_{\max} \tanh\left((\hat B \bv)(j_1,\ldots,j_d)/\tilde u_{\max}\right)
	\label{eq:conconst}
	\end{equation}
	in the TT format using again the TT-Cross method.
	Since the TT approximation may over- or undershoot the exact limits by the relative approximation error $\varepsilon\le \delta$,
	we seek a slightly tighter bound $\tilde u_{\max} = (1-\delta)u_{\max}$.
	
	\subsection{A shifted iterative tensor algorithm for solving \eqref{eq:discnd}}
	\label{sec:iteration}
        The policy update solves \eqref{eq:discnd} at every iteration by taking the previous iterate of the value tensor $\bv$, constructing the control signal, the stiffness matrix and the right hand side, and finally by solving the linear system on the new value tensor approximation.
        The latter step implies using only iterative methods that can preserve the TT structure of all data.
        One of the most robust techniques used nowadays is the Alternating Linear Scheme (ALS) \cite{holtz-ALS-DMRG-2012}
        and its enhanced version, the Alternating Minimal Energy (AMEn) algorithm \cite{ds-amen-2014}.
        In this section we recall the main idea of those, and also propose a shifted version to make the AMEn method applicable to the stationary HJB PDE.

The ALS is a linear projection method similarly to the Krylov techniques,
but in contrast to the latter it projects the equations onto bases constructed from the TT decomposition of the solution itself.
Notice that the TT decomposition \eqref{eq:tt} is linear with respect the the elements of each particular TT block, e.g. $\bv^{(k)}$.
Given \eqref{eq:tt}, let us define a partial TT decomposition where $\bv^{(k)}$ is replaced by the identity matrix $I\in\mathbb{R}^{n \times n}$,
\begin{equation}
\begin{split}
V_{\neq k}(i_1,\ldots,i_d;~\alpha_{k-1},j_k,\alpha_k) & = \sum_{\substack{\alpha_0,\ldots,\alpha_{k-2},\\\alpha_{k+1},\ldots,\alpha_d}} \bv_{\alpha_0,\alpha_1}^{(1)}(i_1) \cdots \bv^{(k-1)}_{\alpha_{k-2},\alpha_{k-1}}(i_{k-1}) \\
& \cdot I(i_k,j_k) \\
& \cdot \bv^{(k+1)}_{\alpha_k,\alpha_{k+1}}(i_{k+1}) \cdots \bv^{(d)}_{\alpha_{d-1},\alpha_d}(i_d).
\end{split}
 \label{eq:frame}
\end{equation}
Clearly, the original TT decomposition \eqref{eq:tt} can be produced from $V_{\neq k}$ by just multiplying it with the $k$-th TT block.
Specifically, we introduce the vector form
\begin{equation}
\bar v^{(k)}(\alpha_{k-1},i_k,\alpha_k) = \bv^{(k)}_{\alpha_{k-1},\alpha_k}(i_k), \quad \bar v^{(k)} \in \mathbb{R}^{r_{k-1}n r_k},
\label{eq:vectorblock}
\end{equation}
and treat $V_{\neq k}$ as a $n^d \times (r_{k-1} n r_k)$ matrix.
One can check that
$$
\tilde \bv = V_{\neq k} \bar v^{(k)}
$$
holds for any $k=1,\ldots,d.$

Now, assuming that the entire stiffness matrix and the right hand side in \eqref{eq:discnd} are assembled in TT formats \eqref{eq:ttm} and \eqref{eq:tt},
the ALS method iterates over $k=1,\ldots,d$ (hence the name alternating), seeking for one TT block at a time by making the residual orthogonal to $V_{\neq k}$,
\begin{equation}\label{eq:localsys}
 \left(V_{\neq k}^\top A V_{\neq k}\right) \bar v^{(k)} = \left(V_{\neq k}^\top \mathbf{b}\right).
\end{equation}
Notice that the reduced system \eqref{eq:localsys} is of size $r_{k-1} n r_k$, i.e. much smaller than the original system \eqref{eq:discnd}.
Once we solve \eqref{eq:localsys}, we populate the TT block $\bv^{(k)}$ with the elements of $\bar v^{(k)}$ through \eqref{eq:vectorblock}, and use the updated $\bv^{(k)}$ to construct \eqref{eq:frame} in the next step $k\rightarrow k+1$ or $k\rightarrow k-1$.
Efficient practical implementation employs the fact that \eqref{eq:frame} mimics the original TT decomposition \eqref{eq:tt} in the sense that same original indices $i_1,\ldots,i_d$ are separated.
Therefore, given \eqref{eq:ttm}, the reduced matrix $V_{\neq k}^\top A V_{\neq k}$ can be computed block by block in a total of $\mathcal{O}(dn^2r^4)$ operations \cite{holtz-ALS-DMRG-2012}.
Here we assume also that $R = \mathcal{O}(r)$, which is the case for the quadratic HJB equation.
Since the reduced matrix in \eqref{eq:localsys} inherits the product structure of the matrix TT decomposition \eqref{eq:ttm}, we can solve \eqref{eq:localsys} iteratively using a fast matrix-vector product with the same cost of $\mathcal{O}(dn^2 r^4)$.

In the AMEn method~\cite{ds-amen-2014},
the TT blocks are also enriched with auxiliary vectors, such as approximate residuals,
which gives a mechanism for increasing TT ranks and adapting them to the desired accuracy.
Consider an auxiliary TT decomposition approximating the residual, projected onto the first $k-1$ TT blocks of the solution,
        \begin{eqnarray*}
                (\mathbf{b} - A\mathbf{v})(i_1,\ldots,i_d)
                & \approx & \sum_{\alpha_0,\ldots,\beta_d=1}^{r_0,\ldots,\rho_d} \bv^{(1)}_{\alpha_0,\alpha_1}(i_1) \cdots \bv^{(k-1)}_{\alpha_{k-2},\alpha_{k-1}}(i_{k-1}) \\
                & \cdot & \mathbf{z}^{(k)}_{\alpha_{k-1},\beta_k}(i_k) \cdot \mathbf{z}^{(k+1)}_{\beta_k,\beta_{k+1}}(i_{k+1}) \cdots \mathbf{z}^{(d)}_{\beta_{d-1},\beta_d}(i_d).
        \end{eqnarray*}
        After solving \eqref{eq:localsys} and updating \eqref{eq:vectorblock},
        we expand $\bv^{(k)}$ and $\bv^{(k+1)}$ increasing their sizes (the corresponding TT ranks) by $\rho_k$,
        \begin{equation}
        \bv^{(k)}(i_k):=\begin{bmatrix}\bv^{(k)}(i_k) & \mathbf{z}^{(k)}(i_k)\end{bmatrix}, \qquad \bv^{(k+1)}(i_{k+1}):=\begin{bmatrix}\bv^{(k+1)}(i_{k+1}) \\ \b0\end{bmatrix}.
        \label{eq:enrich}
        \end{equation}
        Note that this step does not perturb the whole solution tensor $\tilde\bv$ due to the zero block in $\bv^{(k+1)}$, but the subspace of columns of $V_{\neq k+1}$ in the next step is enriched due to the residual block $\mathbf{z}^{(k)}$.
        To reduce the TT rank, we can truncate $\bv^{(k)}$ using the SVD.

Finally, we can \emph{orthogonalise} TT blocks using QR decompositions \cite{osel-tt-2011} such that
\begin{equation}
\sum_{\alpha_{k-1} i_k}\bv^{(k)}_{\alpha_{k-1},\alpha_k}(i_k) \bv^{(k)}_{\alpha_{k-1},\beta_k}(i_k) = I(\alpha_k,\beta_k)
\label{eq:qr-left}
\end{equation}
or
\begin{equation}
\sum_{i_k \alpha_{k}}\bv^{(k)}_{\alpha_{k-1},\alpha_k}(i_k) \bv^{(k)}_{\beta_{k-1},\alpha_k}(i_k) = I(\alpha_{k-1},\beta_{k-1}).
\label{eq:qr-right}
\end{equation}
This makes the projection matrix $V_{\neq k}$ orthogonal as well,
which improves numerical stability of the algorithm.

If the matrix $A$ was symmetric positive definite then the projected system \eqref{eq:localsys}
could be rigorously related to the energy optimization problem and the nonlinear block Gauss--Seidel method.
In our problem \eqref{eq:discnd} this is not the case: $A$ is non-symmetric and degenerate due to the gradient operator $D$, which annihilates any constant component in the solution.
In this case, the degenerate reduced matrix in \eqref{eq:localsys} can prevent convergence.

However, $A$ is compatible with the right hand side under the condition $\ell(\b0)=u_s(\b0) = 0$.
Moreover, the eigenvalues of $A$ are located in the right half of the complex plane for a suitable choice of the domain size $a$ and polynomial order $n$.
Here we propose a modification of AMEn to
resolve both issues by solving shifted systems, mimicking the implicit Euler time propagation.
We introduce a shift $\mu>0$, and solve
\begin{equation}
 \left(A + \mu I \right) \bv = \mathbf{b} + \mu \check \bv,
 \label{eq:shiftnd}
\end{equation}
where $\check \bv$ is the previous iterate of $\bv$.
In practice, we can even combine this shifted AMEn solver and the policy updates into a single iteration as shown in Alg.~\ref{alg:als}.
\begin{algorithm}[h!]
\caption{Policy update with shifted AMEn linear solver}
\label{alg:als}
\begin{algorithmic}[1]
 \State Choose initial value tensor $\bv$, shift $\mu>0$, stopping threshold $\delta>0$, previous iterate $\check\bv=0$, shift reduction factor $0<q<1$.
 \State \textbf{while} $\|\bv - \check\bv\|_2 > \delta \|\bv\|_2$ \textbf{do} \hfill \Comment{Policy iteration}
 \State\quad Set $\check \bv = \bv$ and (optionally) $\mu:=\mu q$.
 \State\quad Compute the control $\mathbf{u}$ using \eqref{eq:Bcontrol} and (optionally) \eqref{eq:conconst}.
 \State\quad  Construct $A[\mathbf{u}]$ and $\mathbf{b}[\mathbf{u}]$ for \eqref{eq:discnd} using Alg.~\ref{alg:ttcross} and \eqref{eq:ttm}--\eqref{eq:ttm_fpq}.
 \State\quad Orthogonalise TT blocks of $\bv$ s.t. \eqref{eq:qr-right} holds for $k=2,\ldots,d$.
 \State\quad  \textbf{for} $k=1,\ldots,d$ \textbf{do} \hfill \Comment{AMEn algorithm}
   \State\quad\quad Assemble and solve $\left(V_{\neq k}^\top A V_{\neq k} + \mu I\right) \bar v^{(k)} = V_{\neq k}^\top \mathbf{b} + \mu V_{\neq k}^\top \check \bv$ using \eqref{eq:frame}.
   \State\quad\quad Update $\bv^{(k)}$ via \eqref{eq:vectorblock} and  truncate $r_k$ using SVD up to accuracy $\delta$.
   \State\quad\quad Enrich $\bv^{(k)},\bv^{(k+1)}$ using \eqref{eq:enrich} and orthogonalise $\bv^{(k)}$ s.t. \eqref{eq:qr-left} holds.
 \State\quad \textbf{end for}
 \State \textbf{end while}
\end{algorithmic}
\end{algorithm}

If $\mathrm{Re}~\lambda (A) \ge 0$, then the spectral radius of the transition matrix $\mu (A+\mu I)^{-1}$ is less than $1$ for any $\mu>0$.
On the other hand, if $\mu > -\bv^\top A \bv$ for any $\bv: \|\bv\|_2=1$,
then the reduced matrix $V_{\neq k}^\top A V_{\neq k} + \mu I$ (remember that $V_{\neq k}^\top V_{\neq k}=I$) is invertible.
This gives a freedom to choose $\mu$ such that the method remains stable and converges fast enough.
In practice we need to ensure $\mu > -\bv^\top A \bv$ only for those $\bv$ that belong to $\mathrm{span}~V_{\neq k}$.
It turns out that as the solution converges, we can decrease $\mu$ geometrically (in particular, we multiply it by a factor $q=0.98$ in each iteration), which ensures faster convergence near the end of the process.

	\section{Optimal feedback control of nonlinear PDEs}\label{spdes} 
\begin{figure}[h!]
	\centering
	\caption{Flowchart illustrating the proposed approach to optimal feedback control of PDEs. The optimal feedback synthesis for the abstract PDE dynamics leads to an infinite-dimensional HJB equation. By using a pseudospectral discretization in space, we obtain a finite-dimensional approximation of the original PDE dynamics, for which the optimal feedback synthesis is realized by solving a finite-dimensional HJB PDE. At this level, we apply our tensor decomposition-based HJB solver, which ultimately leads to a finite-dimensional optimal feedback law approximating the optimal law for the abstract dynamics.}
	\label{fig:diagram}
	\includegraphics[width=\linewidth]{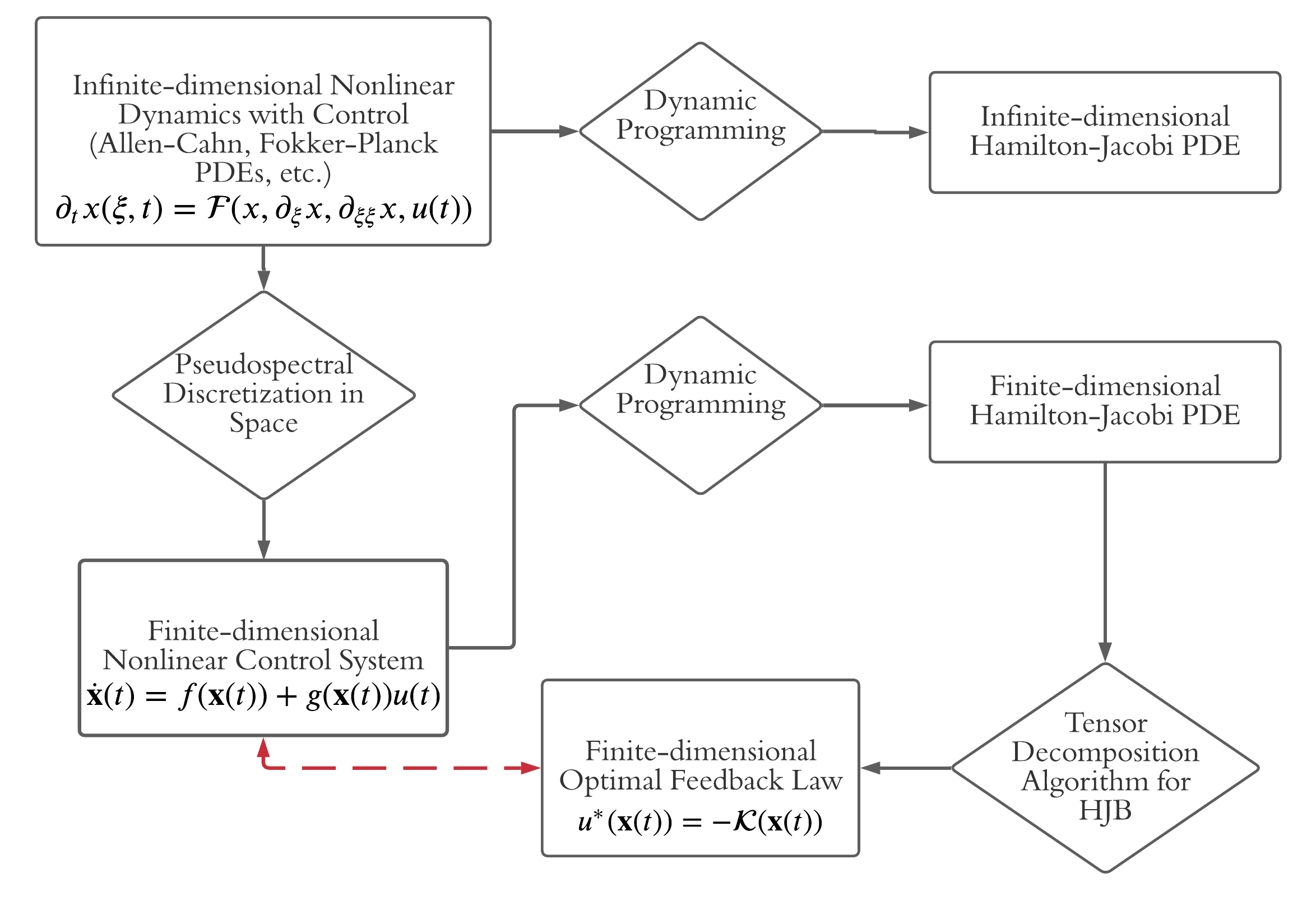}
\end{figure}

High-dimensional nonlinear control systems naturally arise when the dynamics of the system are governed by partial differential equations. From a dynamical perspective, PDEs correspond to abstract, infinite-dimensional systems and therefore the HJB synthesis is understood over an infinite-dimensional state space. Computationally, the treatment is based on the so-called method of lines \cite{mol}. Given an evolutionary PDE, we discretize the space dependence either by finite differences/elements or spectral methods, leading to a large-scale dynamical system with as many state variables as the space discretization dictates. We perform the HJB synthesis over this finite but high-dimensional system, as summarized in Fig.\ref{fig:diagram}. The accuracy of such a representation and its implications over the control design vary depending on the class of PDEs under consideration. Strongly dissipative PDEs can be accurately represented with few degrees of freedom in space, while convection-dominated PDEs might require a much more complex state space representation. Therefore, the performance study of our framework with respect to the dimension parameter is central to our analysis. It benefits from the fact that, unlike the nonlinear ODE world, the taxonomy of physically meaningful nonlinearities in time-dependent PDEs is well-delimited. We focus on nonlinear reaction PDEs where we take the Allen-Cahn equation as a reference model due to its rich equilibrium structure, and to nonlinear convection in the Fokker-Planck equation, where the control action enters as a bilinear term. The semi-discretization in space of these nonlinearities leads to well-structured finite-dimensional realizations \cite{KK18} which allow a systematic analysis of the scaling of our methodology with respect to the dimension.
	A Matlab implementation of Algorithm~\ref{alg:als} and the numerical examples benchmarked below are available at \url{https://github.com/dolgov/TT-HJB}.
	All tests were carried out in Matlab 2017b on one core of an Intel Xeon E5-2640 v4 CPU at 2.4 GHz.

	\subsection{The Allen-Cahn equation}
	We consider the following nonlinear diffusion-reaction equation \cite{KK18}:
	\begin{equation}
	\begin{split}
	\partial_t x(\xi,t) & = \sigma \partial_{\xi\xi} x + x(1-x^2) + \chi_{\omega}(\xi) u(t)\,,
	\end{split}
	\label{eq:nw}
	\end{equation}
	in  $[-1,1] \times \mathbb{R}^+$ with Neumann boundary conditions. We set $\sigma=0.2$, and the scalar control signal $u(t)$ acts through the indicator $\chi_{\omega}$ of the subdomain $\omega=[-0.5, 0.2]$. This equation has $x=0$ as unstable equilibrium and $x=\pm 1$ as stable equilibria. The cost is given by
	\begin{equation}
	\mathcal{J}(u,x) = \int_{0}^{\infty} \|x(\xi,t)\|_{L_2(-1,1)}^2 + \gamma u(t)^2 dt\,.
	\label{eq:cost_unc}
	\end{equation}
	where for concreteness we take $\gamma=0.1$. Thus, the control objective consists in stabilizing the unstable equilibrium. \eqref{eq:nw} is discretized by Chebyshev pseudospectral collocation method \cite{trefethen-spectral-2000} using $d$ points $\xi_k = -\cos(\pi k /(d+1))$, $k=1,\ldots,d$.
	The discrete state is collected into a vector of nodal values $X(t) = \left(X_1(t),\ldots,X_d(t)\right)$
	where $X_k(t) \approx x(\xi_k,t)$, leading to a  $d$-dimensional nonlinear ODE
	\begin{equation}
	\frac{dX}{dt} = A X + X \odot (1-X\odot X) + Bu(t),
	\label{eq:nw_discr}
	\end{equation}
	where ``$\odot$'' is the coordinatewise Hadamard product, $A$ is the pseudospectral differentiation matrix corresponding to the Laplace operator, and $B$ is a vector corresponding to the pseudospectral discretisation of the indicator function $\chi_{\omega}(\xi)$.
	The HJB equation solver is applied directly to \eqref{eq:nw_discr}, restricting the domain of the value function to $(-3,3)^d$,
	sufficient to accommodate typical initial states.
	We compare our design against the linear-quadratic regulator LQR feedback law \cite[Chapter 8]{LQR} computed for the dynamics in \eqref{eq:nw_discr} linearised around the origin, $A_{L} = A + I$.
	
	\subsubsection*{Algorithmic performance for the one-dimensional Allen-Cahn equation}
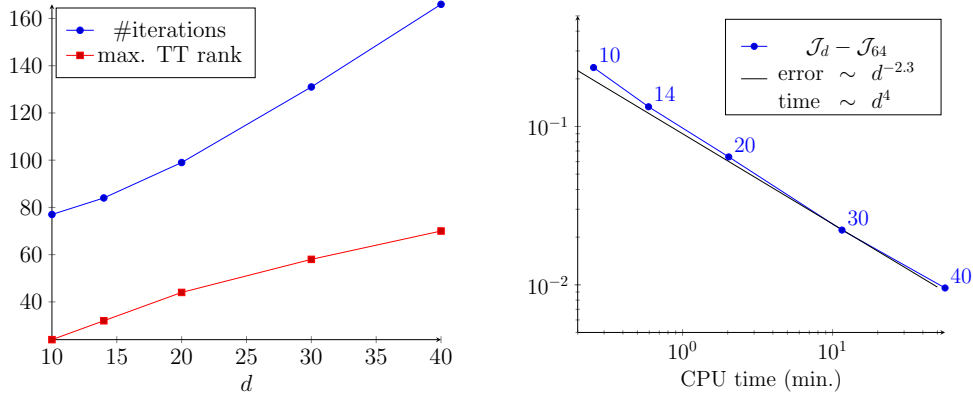
\begin{figure}[h!]
		\centering
		\caption{Allen-Cahn problem \eqref{eq:nw}. Left: numbers of policy iterations and maximal TT ranks. Right: differences in total running cost and CPU times for different spatial dimensions~$d$. Numbers above points in the right plot denote~$d$.  Other settings $n=5$ and $\delta=10^{-3}$.}
		\label{fig:nw-d}
			\resizebox{0.47\linewidth}{!}{%
				\begin{tikzpicture}
				\begin{axis}[%
				xmode=normal,
				ymode=normal,
				xlabel={$d$},
				legend style={inner sep=6pt,at={(0.01,0.99)},anchor=north west},
				label style={font=\Large},
				tick label style={font=\Large}
				]
				\addplot+[] coordinates{
					(10, 77 )
					(14, 84 )
					(20, 99 )
					(30, 131)
					(40, 166)
				}; \addlegendentry{\Large \#iterations};
				\addplot+[] coordinates{
					(10, 24)
					(14, 32)
					(20, 44)
					(30, 58)
					(40, 70)
				}; \addlegendentry{\Large max. TT rank};
				\end{axis}
				\end{tikzpicture}
		}\hfill\resizebox{0.47\linewidth}{!}{%
				\begin{tikzpicture}
				\begin{axis}[%
				xmode=log,
				ymode=log,
				ymin=5e-3,ymax=5e-1,
				xlabel={CPU time (min.)},
				legend style={inner sep=10pt,at={(0.99,0.99)},anchor=north east},
				nodes near coords,point meta=explicit symbolic,
				every node near coord/.append style={anchor=south west},
				label style={font=\Large},
				tick label style={font=\Large}
				]
				
				\addplot+[] coordinates{
					(15.309717 /60, 1.772745-1.537047)[\Large 10]
					(35.609753 /60, 1.670270-1.537047)[\Large 14]
					(121.982704/60, 1.601298-1.537047)[\Large 20]
					(693.843014/60, 1.559215-1.537047)[\Large 30]
					(3374.85754/60, 1.546595-1.537047)[\Large 40]
				}; \addlegendentry{\Large $\mathcal{J}_d - \mathcal{J}_{64}$};
				\addplot+[no marks,domain=0.2:50,black] {0.09*x^(-0.57)}; \addlegendentry{\Large $\begin{array}{lcl}\mbox{error}&\sim &d^{-2.3} \\ \mbox{time}&\sim &d^4\end{array}$};
				\end{axis}
				\end{tikzpicture}
			}
	\end{figure}

	We first investigate the performance of Algorithm~\ref{alg:als} with respect to the dimension of the dynamical system.
	We fix $n=5$, $\delta=10^{-3}$, and the initial shift in Alg.~\ref{alg:als} $\mu=50$.
	In Fig.~\ref{fig:nw-d} we can observe that the maximal TT rank grows linearly with the dimension. The $\mathcal{O}(dn^2r^4)$ complexity of the AMEn method leads to a total cost bound in order of $d^5$. However, the
	effective cost is closer to $\mathcal{O}(d^4)$ (Fig.~\ref{fig:nw-d}, right), which can be attributed to a non-uniform distribution of TT ranks along the decomposition.
	This is a significant reduction compared to the exponential cost of the full Cartesian ansatz $n^d$.
	However, the method can become slow in very high dimensions,
	mainly due to the increase in the number of policy iterations as shown in Fig.~\ref{fig:nw-d} (left), resulting from a larger condition number of the linearised system.

        \begin{figure}[h!]
                \centering
                \caption{Allen-Cahn problem \eqref{eq:nw} with $d=40$, $n=5$ and $\delta=10^{-3}$. Left: time evolution of running costs. Right: control signals. The uncontrolled state converges to $X=1$ with an infinite cost.}
                \label{fig:nw-ut}
		\resizebox{0.47\linewidth}{!}{%
				\begin{tikzpicture}
				\begin{axis}[%
				xmode=normal,
				ymode=log,
				xlabel={$t$},
				xmin=0,xmax=3.2,
				legend style={at={(0.01,0.01)},anchor=south west},
				label style={font=\Large},
				tick label style={font=\Large}
				]
				\addplot+[no marks,line width=1.0pt] table[header=true,x=t,y=chjb]{nw_cost_d40.dat} node[pos=0.3,anchor=east,inner sep=1pt]{\Large $\begin{array}{l}\mbox{HJB}_{40}\\\mathcal{J}=1.55\end{array}$};
				\addplot+[no marks,line width=1.5pt,dashed] table[header=true,x=t,y=c40]{nw_cost_hjb_interp14.dat} node[pos=0.25,anchor=south west,inner sep=1pt]{\Large  $\begin{array}{l}\mbox{HJB}_{14}\\\mathcal{J}=1.73\end{array}$};
				\addplot+[no marks,line width=2.5pt] table[header=true,x=t,y=clqr]{nw_cost_d40.dat} node[pos=0.47,anchor=north,inner sep=1pt]{\Large  $\begin{array}{l}\mbox{LQR}\\\mathcal{J}=40.2\end{array}$};
				\addplot+[no marks,line width=3.0pt,opacity=0.5] table[header=true,x=t,y=cunc]{nw_cost_d40.dat} node[pos=0.3,anchor=south west,inner sep=1pt]{\Large  $\begin{array}{l}\mbox{UNC}\\\mathcal{J}=\infty\end{array}$};
				\end{axis}
				\end{tikzpicture}
		}
		\hfill\resizebox{0.47\linewidth}{!}{%
				\begin{tikzpicture}
				\begin{axis}[%
				xmode=normal,
				ymode=normal,
				xlabel={$t$},
				xmin=0,xmax=3.2,
				ymin=-10,ymax=5,
				legend style={at={(0.01,0.01)},anchor=south west},
				label style={font=\Large},
				tick label style={font=\Large}
				]
				\addplot+[no marks,line width=1.0pt] table[header=true,x=t,y=uhjb]{nw_cost_d40.dat} node[pos=0.15,anchor=north,inner sep=10pt]{\Large $\mbox{HJB}_{40}$};
				\addplot+[no marks,line width=1.5pt,dashed] table[header=true,x=t,y=u40]{nw_cost_hjb_interp14.dat} node[pos=0.08,anchor=south,inner sep=10pt]{\Large $\mbox{HJB}_{14}$};
				\addplot+[no marks,line width=2.5pt] table[header=true,x=t,y=ulqr]{nw_cost_d40.dat} node[pos=0.7,anchor=west]{\Large LQR};
				\end{axis}
				\end{tikzpicture}
			}
	\end{figure}
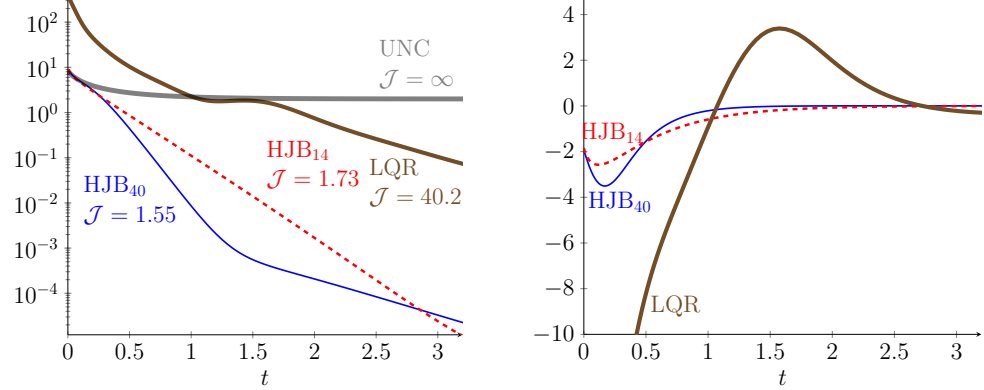

        A possible remedy is to construct the value function for a lower-dimensional discretisation of the PDE,
        and interpolate the state of a system with finer discretisation onto this lower-dimensional spatial grid.
        The resulting error, proportional to the discretisation error of the lower-dimensional grid,
                might be still much smaller than the error resulting from e.g. the linearisation of the system.
	Figure~\ref{fig:nw-ut} (left) shows the running costs $\mathcal{J}=\|x(\xi,t)\|^2 + \gamma u(t)^2$ for the HJB control with system dimension $40$ with $n=5$ and $\delta=10^{-3}$ ($\mbox{HJB}_{40}$), an interpolated HJB control from dimension $14$ and the same $n=5,\delta=10^{-3}$ ($\mbox{HJB}_{14}$),
	an LQR feedback, and finally the cost of the UNControlled system (both with $d=40$),
	with the initial condition $x(\xi,0) = 2+\cos(2\pi\xi)\cos(\pi\xi)$.
	Fig.~\ref{fig:nw-ut} (right) shows the corresponding control signals.
	Since the linearized system is unstable, the LQR acts very aggressively during the transient phase.
	The HJB synthesis is able to detect the stabilizing  effect of the nonlinearity and produce a control at much lower cost. We observe differences in the control signals and total costs for the HJB laws depending on the dimension of the dynamical system, justifying the need for accurate high-dimensional HJB solvers.

        \begin{figure}[h!]
                \centering
                \caption{Allen-Cahn problem \eqref{eq:nw} with $\delta=10^{-3}$ and $d=14$. Left: numbers of policy iterations and maximal TT ranks. Right: differences in total running cost and CPU times for different univariate polynomial degrees~$n$. Numbers above points in the right plot denote~$n$.}
                \label{fig:nw-n}
                \noindent\resizebox{0.47\linewidth}{!}{%
                        \begin{tikzpicture}
                        \begin{axis}[%
                        xmode=normal,
                        ymode=normal,
                        xlabel={$n$},
                        xmin=3,xmax=7,
                        xtick={3,4,5,6,7},
                        legend style={at={(0.01,0.99)},anchor=north west},
                        ]
                        \addplot+[] coordinates{
                                (3, 30  )
                                (4, 58  )
                                (5, 61  )
                                (6, 118 )
                                (7, 150 )
                        }; \addlegendentry{\#iterations};
                        \addplot+[] coordinates{
                                (3, 13 )
                                (4, 27 )
                                (5, 33 )
                                (6, 56 )
                                (7, 90 )
                        }; \addlegendentry{max. TT rank};
                        \end{axis}
                        \end{tikzpicture}
                }
                \resizebox{0.47\linewidth}{!}{%
                        \begin{tikzpicture}
                        \begin{axis}[%
                        xmode=log,
                        ymode=log,
                        ymin=5e-3,ymax=2e1,
                        xlabel={CPU time (min.)},
                        legend style={at={(0.99,0.99)},anchor=north east},
                        nodes near coords,point meta=explicit symbolic,
                        every node near coord/.append style={anchor=south west},
                        ]

                        \addplot+[] coordinates{
                                (3.595657  /60, 9.763848-1.544810)[3]
                                (13.370541 /60, 2.389479-1.544810)[4]
                                (24.752617 /60, 1.669412-1.544810)[5]
                                (243.833027/60, 1.564735-1.544810)[6]
                                (1967.37   /60, 1.544810-1.533555)[7]
                        }; \addlegendentry{$\mathcal{J}_n - \mathcal{J}_8$};

                        \end{axis}
                        \end{tikzpicture}
                }
        \end{figure}
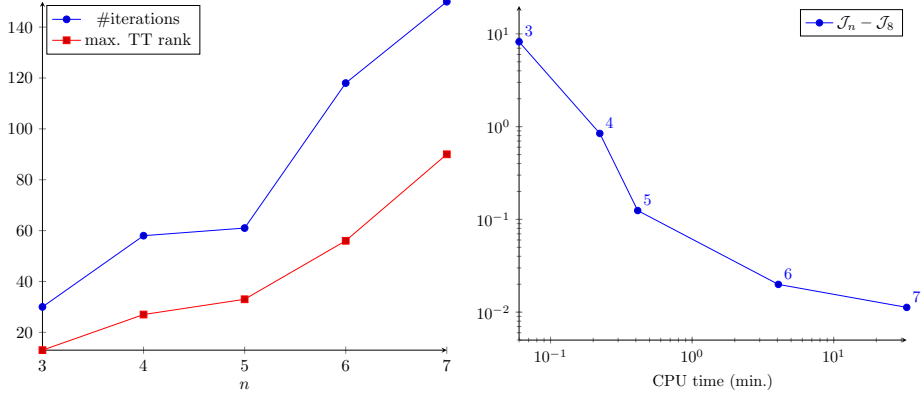

        \begin{figure}[h!]
                \centering
                \caption{Allen-Cahn problem \eqref{eq:nw} with $n=5$ and $d=14$. Left: numbers of policy iterations and maximal TT ranks. Right: differences in total running cost and CPU times for different TT approximation thresholds~$\delta$. Numbers above points in the right plot denote~$\delta$.}
                \label{fig:nw-delta}
                \noindent\resizebox{0.47\linewidth}{!}{%
                        \begin{tikzpicture}
                        \begin{axis}[%
                        xmode=log,
                        ymode=normal,
                        xlabel={$\delta$},
                        legend style={at={(0.99,0.99)},anchor=north east},
                        ]
                        \addplot+[] coordinates{
                                (1e-2, 40)
                                (3e-3, 52)
                                (1e-3, 61)
                                (3e-4, 72)
                                (1e-4, 82)
                        }; \addlegendentry{\#iterations};
                        \addplot+[] coordinates{
                                (1e-2, 18)
                                (3e-3, 25)
                                (1e-3, 33)
                                (3e-4, 44)
                                (1e-4, 57)
                        }; \addlegendentry{max. TT rank};
                        \end{axis}
                        \end{tikzpicture}
                }
                \resizebox{0.47\linewidth}{!}{%
                        \begin{tikzpicture}
                        \begin{axis}[%
                        xmode=log,
                        ymode=log,
                        xmin=0.1,xmax=2,
                        ymin=1e-4,ymax=1e-1,
                        xlabel={CPU time (min.)},
                        legend style={at={(0.99,0.99)},anchor=north east},
                        nodes near coords,point meta=explicit symbolic,
                        every node near coord/.append style={anchor=south west},
                        ]

                        \addplot+[] coordinates{
                                (8.971272 /60, 1.721134-1.668442)[$1\cdot 10^{-2}$]
                                (15.916878/60, 1.677024-1.668442)[$3\cdot 10^{-3}$]
                                (24.752617/60, 1.669412-1.668442)[$1\cdot 10^{-3}$]
                                (58.557380/60, 1.668442-1.668198)[$3\cdot 10^{-4}$]
                        }; \addlegendentry{$\mathcal{J}_{\delta} - \mathcal{J}_{10^{-4}}$};
                        \addplot+[no marks,domain=0.1:2,black] {4.63e-5*x^(-3)}; \addlegendentry{$\mbox{Time}^{-3}$};
                        \end{axis}
                        \end{tikzpicture}
                }
        \end{figure}
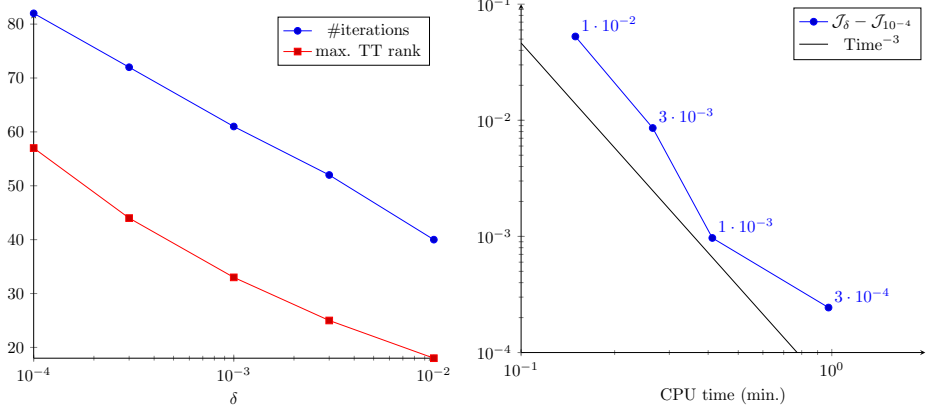

        Now we analyse the performance of the TT-HJB scheme depending on the number of Legendre polynomials $n$ in each variable (Fig. \ref{fig:nw-n}) and the stopping threshold $\delta$ in Alg.~\ref{alg:als} (Fig. \ref{fig:nw-delta}).
        Due to the nonlinearity in \eqref{eq:nw_discr}, the value function is significantly far from a quadratic polynomial, which is reflected in Fig. \ref{fig:nw-n} by the linear growth of TT ranks with $n$, and a relatively slow algebraic decay of the error.
        Nevertheless, even an order-4 approximation can give a substantially better control signal than the LQR approximation, see Fig.~\ref{fig:nw-ut}.
        From Fig.~\ref{fig:nw-delta} we see that the number of iterations and the TT ranks depend logarithmically on the TT approximation error, which is a more optimistic result than that predicted by Thm.~\ref{thm:rank}, although the problem is nonlinear.
        The errors in the total cost start higher than the threshold $\delta$, but eventually the two error indicators are of the same order.

	\subsubsection*{Allen-Cahn problem with control constraints}
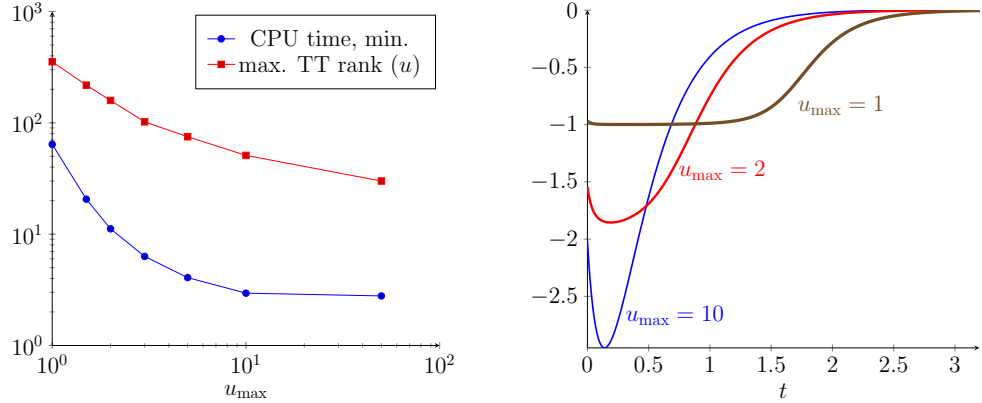
\begin{figure}[h!]
		\centering
		\caption{Allen-Cahn problem with control constraints. Left: CPU times and TT ranks. Right: control signals for different control constraints $-u_{\max} \le u \le u_{\max}$. We set $d=20$, $\delta=10^{-3}$ and $n=5$.}
		\label{fig:nw-umax}
		\resizebox{0.47\linewidth}{!}{%
				\begin{tikzpicture}
				\begin{axis}[%
				xmode=log,
				ymode=log,
				xlabel={$u_{\max}$},
				ymin=1,ymax=1e3,
				xmin=1,xmax=100,
				legend style={inner sep=7pt,at={(0.99,0.99)},anchor=north east},
				label style={font=\Large},
				tick label style={font=\Large}
				]
				\addplot+[] coordinates{
					(50, 167.236885/60)
					(10, 177.062251/60)
					(5 , 244.297876/60)
					(3 , 378.799143/60)
					(2 , 669.878672/60)
					(1.5,1237.13782/60)
					(1 , 3847.410  /60)
				}; \addlegendentry{\Large CPU time, min.};
				\addplot+[] coordinates{
					(50, 30 )
					(10, 51 )
					(5 , 75 )
					(3 , 102)
					(2 , 159)
					(1.5,218)
					(1 , 354)
				}; \addlegendentry{\Large max. TT rank ($u$)};
				
				\end{axis}
				\end{tikzpicture}
			}
		\hfill\resizebox{0.47\linewidth}{!}{%
				
				\begin{tikzpicture}
				\begin{axis}[%
				xmode=normal,
				ymode=normal,
				xmin=0,xmax=3.2,
				xlabel={$t$},
				legend style={at={(0.99,0.99)},anchor=north east},
				label style={font=\Large},
				tick label style={font=\Large}
				]
				
				\addplot+[no marks,line width=1.0pt] table[header=true,x=t,y=uhjb]{nw_umax10_d20.dat} node[pos=0.1,anchor=west]{\Large $u_{\max}=10$};
				\addplot+[no marks,line width=1.5pt] table[header=true,x=t,y=uhjb]{nw_umax2_d20.dat} node[pos=0.10,anchor=west]{\Large $u_{\max}=2$};
				\addplot+[no marks,line width=2.0pt] table[header=true,x=t,y=uhjb]{nw_umax1_d20.dat} node[pos=0.17,anchor=north west]{\Large $u_{\max}=1$};
				\end{axis}
				\end{tikzpicture}
			}
	\end{figure}

	As discussed in Section \ref{shjb}, the proposed framework allows to enforce control constraints through a suitable choice of the control penalties in \eqref{optc2}. Figure~\ref{fig:nw-umax} (left) shows the total CPU times and TT ranks of the constrained feedback law.
	Figure~\ref{fig:nw-umax} (right) presents the control signals for three bound parameters. As $\mathcal{P}(x)$ becomes steeper for more severe control constraints, the TT ranks increase leading to longer computing times.
	Nevertheless, Alg.~\ref{alg:als} remains effective for a wide range of constraints, adjusting the value function accordingly.

	\subsubsection*{Allen-Cahn problem with 2-dimensional space}
	We study the extension of the problem \eqref{eq:nw} to two spatial dimensions with the state depending on two space coordinates, $x(\boldsymbol\xi,t) = x(\xi_1,\xi_2,t)$.
	We replace the second derivative with the Laplace operator, and the control is applied on the domain $\omega=[-0.5,0.2]^2$.
	We use the Cartesian product of the same Chebyshev grids in each direction, and similar homogeneous Neumann conditions on the boundary of $\Omega=[-1,1]^2$. The CPU times and TT ranks are shown in Fig.~\ref{fig:nw2-d} (left).
	Although Theorem~\ref{thm:rank} is not immediately applicable to a nonlinear system,
	we still observe a linear growth of the TT ranks with the number of Chebyshev points in each direction.
	The values of the TT ranks are larger than those in the one-dimensional case,
	leading to increased computing times.
	However, the performance of the high-dimensional HJB controller is satisfactory.
	Figures~\ref{fig:nw2-d} (right) and \ref{fig:nw2-ut} show the response of the system with an initial state $x(\boldsymbol\xi,0) = 2+\cos(2\pi\xi_1)\cos(\pi\xi_2).$ We can see again that the HJB-controlled state is stabilized while the LQR synthesis fails.
	
	\begin{figure}[h!]
		\centering
		\caption{Two-dimensional Allen-Cahn control problem. Left: CPU times and TT ranks. Right: HJB control signals.}
		\label{fig:nw2-d}
		\resizebox{0.47\linewidth}{!}{%
				\begin{tikzpicture}
				\begin{axis}[%
				xmode=normal,
				ymode=log,
				xlabel={$d$},
				legend style={inner sep=6pt,at={(0.01,0.99)},anchor=north west},
				label style={font=\Large},
				tick label style={font=\Large}
				]
				\addplot+[] coordinates{
					(7*7 ,  5684.846/60)
					(8*8 ,  21375.38/60)
					(9*9 ,  56203.74/60)
					(10*10, 167236.1/60)
					(11*11, 458381.7/60)
				}; \addlegendentry{\Large CPU time, min.};
				\addplot+[] coordinates{
					(7*7 , 99 )
					(8*8 , 136)
					(9*9 , 163)
					(10*10, 210)
					(11*11, 230)
				}; \addlegendentry{\Large max. TT rank};
				\addplot+[no marks,dashed,black,domain=50:120] {33*sqrt(x)-131}; \addlegendentry{\Large $\mathcal{O}(\sqrt{d})$};
				\end{axis}
				\end{tikzpicture}
			}
		\hfill\resizebox{0.47\linewidth}{!}{%
				\begin{tikzpicture}
				\begin{axis}[%
				xmode=normal,
				ymode=normal,
				xmin=0,xmax=3.2,
				xlabel={$t$},
				legend style={at={(0.99,0.99)},anchor=north east},
				label style={font=\Large},
				tick label style={font=\Large}
				]
				
				\addplot+[no marks,line width=1.0pt] table[header=true,x=t,y=uhjb]{nw2_cost_d9.dat} node[pos=0.32,anchor=west]{\Large $d=81$};
				\addplot+[no marks,line width=1.0pt] table[header=true,x=t,y=uhjb]{nw2_cost_d11.dat} node[pos=0.18,anchor=west]{\Large $d=121$};
				\end{axis}
				\end{tikzpicture}
			}
	\end{figure}
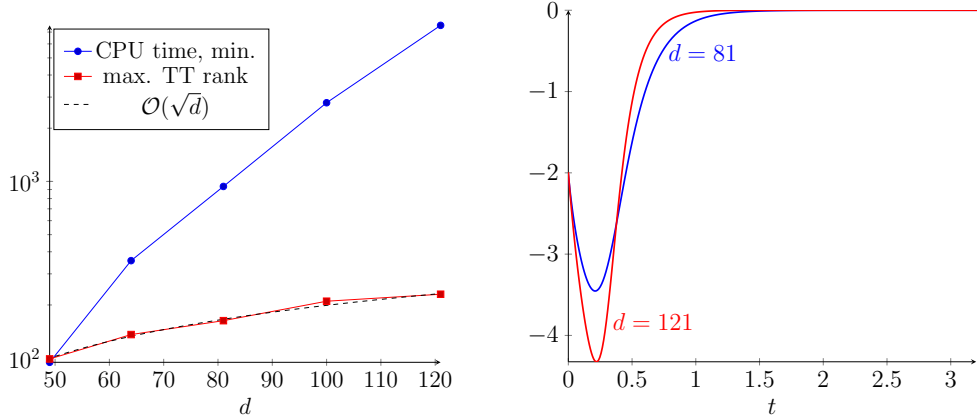
        \begin{figure}[h!]
                \centering
                \caption{Two-dimensional Allen-Cahn control problem with $d=121$. Left: state snapshots at $t=0.6$ for HJB and LQR control laws. Right: total running costs.}
                \label{fig:nw2-ut}
		\includegraphics[width=0.47\linewidth]{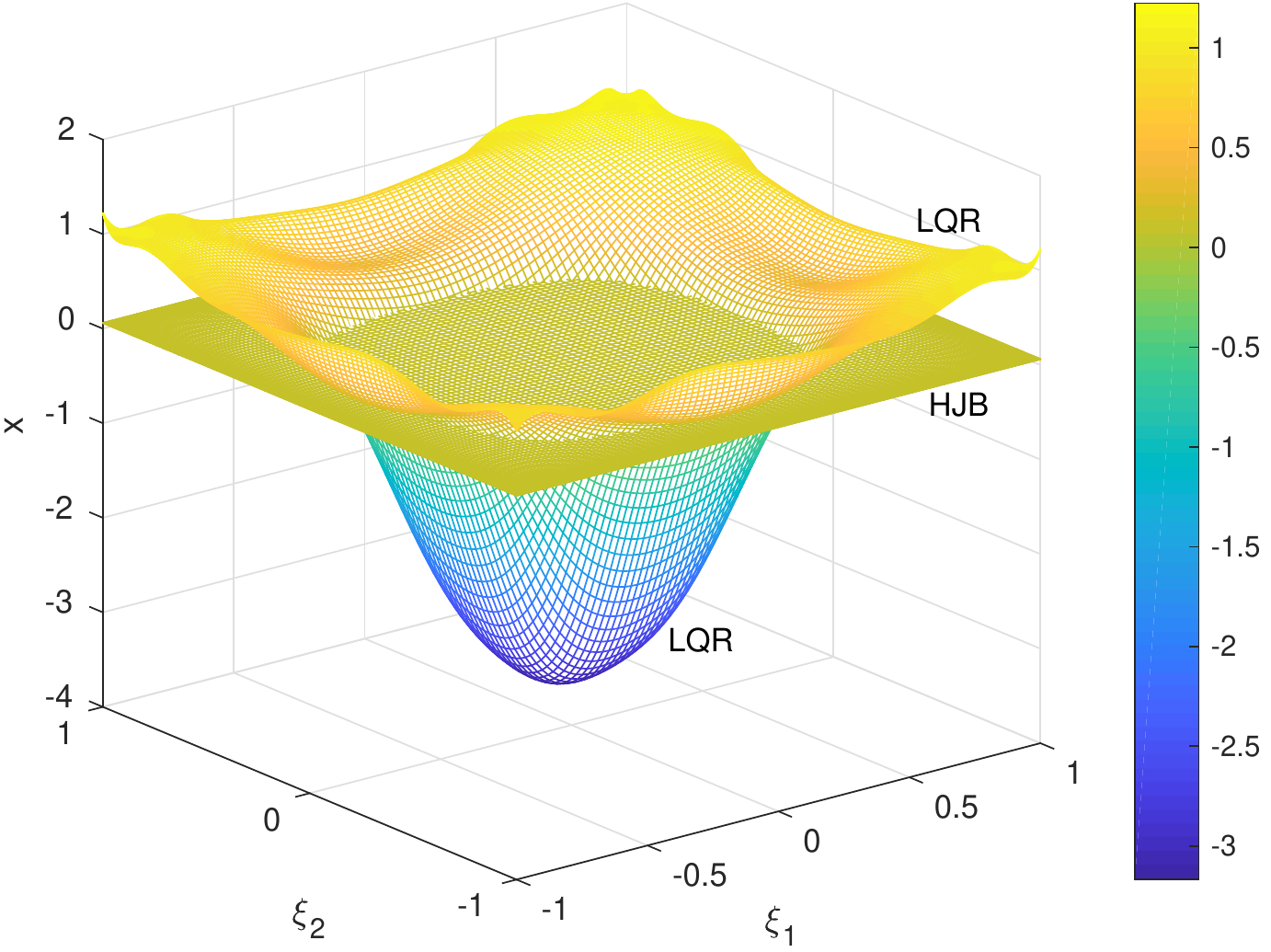}
		\hfill\resizebox{0.47\linewidth}{!}{%
				\begin{tikzpicture}
				\begin{axis}[%
				xmode=normal,
				ymode=log,
				xmin=0,xmax=3.2,
				xlabel={$t$},
				legend style={at={(0.99,0.99)},anchor=north east},
				label style={font=\Large},
				tick label style={font=\Large}
				]
				\addplot+[no marks,line width=1.0pt] table[header=true,x=t,y=chjb]{nw2_cost_d11_new.dat} node[pos=0.32,anchor=west]{\Large HJB};
				\addplot+[no marks,line width=1.0pt] table[header=true,x=t,y=clqr]{nw2_cost_d11_new.dat} node[pos=0.18,anchor=west]{\Large LQR};
				\end{axis}
				\end{tikzpicture}
			}
	\end{figure}
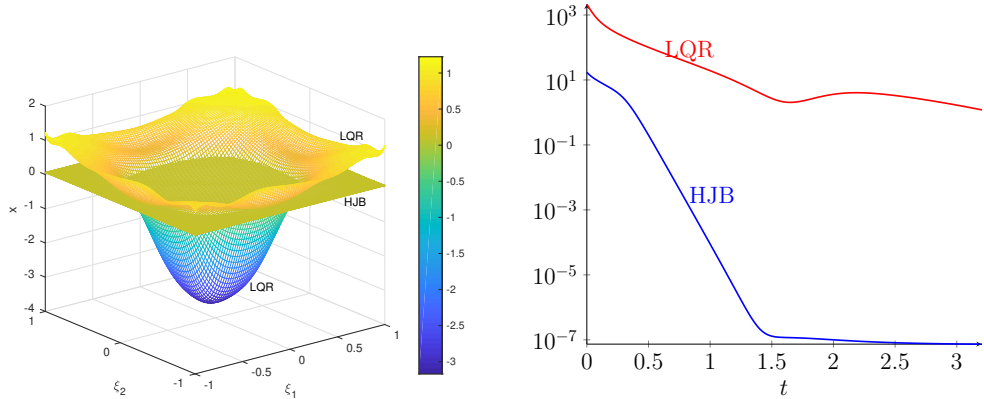

	\subsection{The Fokker-Planck equation}\label{sec:FPnum}
	We compute optimal feedback regulators for the stabilised bilinearly controlled Fokker-Planck equation
	\begin{equation}
	\begin{split}
	\partial_t x(\xi,t) & =  \nu \partial_{\xi\xi} x + \partial_{\xi} (x \partial_{\xi} G) + u \partial_{\xi} (x \partial_{\xi} H), \quad \xi \in \Omega,\\
	0 & = \left[\partial_{\xi} x + x \partial_{\xi} (G+uH)\right]|_{\xi \in \partial \Omega},
	\end{split}
	\label{eq:fpe}
	\end{equation}
	where the computational domain will be set $\Omega=(-6,6).$ This equation models the density of particles,
	controlled with laser-induced electric force with potential $G(x)+u(t) H(x)$ \cite{HarST13}, with $G$ the ground and $H$ the control potential. We refer to \cite{AB18} for a recent survey on further illustrations on the importance of pdf-based method for the control of the Fokker-Planck equation, including the control agent-based dynamics \cite{ACFK}.
	This system  has $0$ as an eigenvalue with associated eigenstate $x_\infty= \exp(- (\log \nu + \frac{G}{\nu}))$, see eg. \cite{bkp-hjb-fp-2018}. Henceforth the eigenstates are considered as normalized in $L^2(\Omega)$.
	It  is known that $x_\infty$ is stable, but the convergence to this steady state,  which is given by the second eigenvalue and  depends on $\nu$ and $G$,
	can be extremely slow, see for instance \cite[pg 251]{MatS81}. Thus, to speed up convergence  in the transient phase, control is of importance.  To obtain a suitable stabilization problem we introduce the shifted state $y=x - x_\infty$. It satisfies
	\begin{equation}
	\begin{split}
	\partial_t y(\xi,t) & =  \nu \partial_{\xi\xi} y + \partial_{\xi} (y \partial_{\xi} G) + u \partial_{\xi} (y \partial_{\xi} H)+ u \partial_{\xi} (x_\infty \partial_{\xi} H), \quad \xi \in \Omega,\\
	0 & = \left[\partial_{\xi} y + y \partial_{\xi} G +u(y+x_\infty) \partial_\xi H \right]|_{\xi \in \partial \Omega}.
	\end{split}
	\label{eq:sfpe}
	\end{equation}
	The control objective consists now in driving $y$ to zero.
	To compute the controller  we further  introduce  a positive, i.e. destabilising,  shift by adding $\sigma y$ to the right hand side of \eqref{eq:sfpe}.
	If  this controller is applied to the unshifted equation it accelerates convergence of $y$ to 0 and hence the convergence of $x$ to $x_\infty$.
	
	Considering the variational form of \eqref{eq:sfpe} one observes that the control will not have an effect on a subspace of co-dimension one. For this reason we introduce $Y_{\mathcal{P}}=\{v\in L^2(\Omega): \int_\Omega v\,d\xi =0\}$, and denote by $ {\mathcal{P}} \in {\mathcal{L}}(L^2(\Omega),Y_{\mathcal{P}})$  the projection onto $Y_{\mathcal{P}}$ along  $x_\infty$, which is given by ${\mathcal{P}}y= y-(\int_\Omega y\, d\xi) x_\infty$.  Subsequently we apply ${\mathcal{P}}$ to  \eqref{eq:sfpe} with initial datum given by ${\mathcal{P}} x(0)$. For the details we refer to \cite{BreKP18}.
	
	The  Fokker-Planck equation \eqref{eq:sfpe} is discretized using a finite difference  scheme with $D$ intervals. To allow for possible further reduction of the dimension a  balanced truncation based model reduction, adapted to bilinear systems \cite{BenD11}, is used, to reduce the system to dimension $d$.
	
\begin{figure}[h!]
		\centering
		\caption{Ground and control potentials in the Fokker-Planck control system.}
		\label{fig:fppot}
		\includegraphics[width=0.7\linewidth]{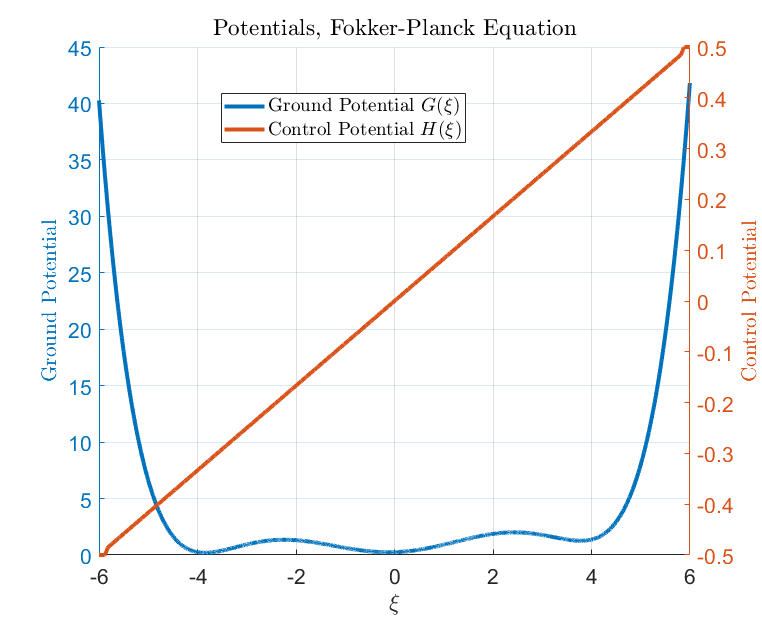}
	\end{figure}

	For the numerical results  we fix $\gamma=10^{-2}, \nu=1,$ $\sigma = 0.2,$ and the potentials $G(\xi)$ and $H(\xi)$ are chosen to reproduce the setting in \cite{bkp-hjb-fp-2018}, as shown in Fig.~\ref{fig:fppot}. That is, the ground potential is set to be
	\begin{equation}
	G(\xi)=\frac{\left((0.5\xi^2-15)\xi^2+119\right)\xi^2+28\xi+50}{200}\,,
	\end{equation}
	whereas $H(\xi)$ is given by
	\begin{equation}
	H(\xi)=\begin{cases}
	-1/2 & \text{if}\quad -6.0\leq\xi\leq-5.9\\
	\xi/12& \text{if}\quad -5.8\leq\xi\leq 5.8\\
	1/2&\text{if}\quad 5.9\leq \xi\leq 6.0
	\end{cases}
	\end{equation}
	with the disjoint intervals united with an Hermite interpolant.

	Since both the original system size $D$, and the reduced dimension $d$, are approximation parameters, we need to set them to appropriate values that deliver a desired accuracy in the model outcomes, such as the total cost. Note also that in contrast to the linear case, the generalized balanced truncation method for bilinear systems does
	not exhibit an a priori error bound \cite{BenD11}.
	
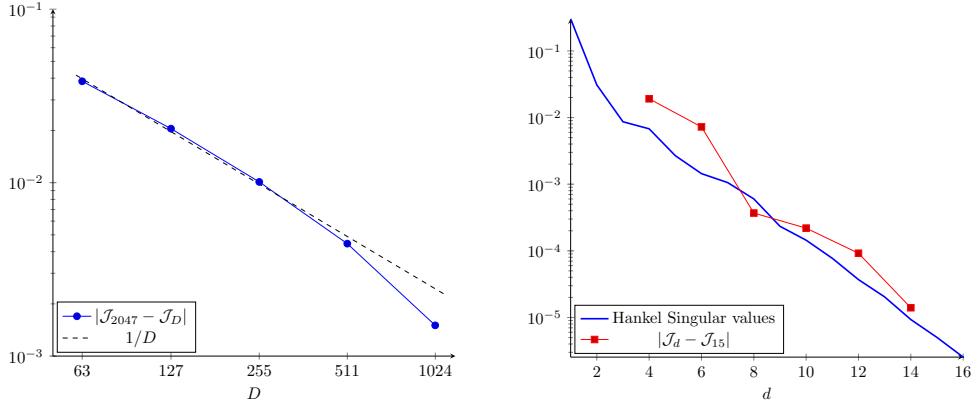
\begin{figure}[h!]
		\centering
		\caption{Errors in the total cost in the Fokker-Planck model with the uniform initial state $x(\xi,0) = \frac{1}{12}$ for different numbers of discretisation points (left) and different dimensions of the reduced model (right).}
		\label{fig:fp-Dd}
		\resizebox{0.47\linewidth}{!}{%
			\begin{tikzpicture}
			\begin{axis}[%
			xmode=log,
			ymode=log,
			xlabel={$D$},
			xmin=50,xmax=1200,
			ymax=1e-1,ymin=1e-3,
			xtick={63,127,255,511,1024},
			xticklabels={63,127,255,511,1024},
			legend style={at={(0.01,0.01)},anchor=south west},
			]
			\addplot+[] coordinates{
				(63  , 0.151559 - 0.113136)
				(127 , 0.151559 - 0.131079)
				(255 , 0.151559 - 0.141463)
				(511 , 0.151559 - 0.147110)
				(1023, 0.151559 - 0.150055)
			}; \addlegendentry{$|\mathcal{J}_{2047}-\mathcal{J}_D|$};
			\addplot+[no marks,black,dashed,domain=60:1100] {2.5*x^(-1.0)}; \addlegendentry{$1/D$};
			\end{axis}
			\end{tikzpicture}
		}
		\hfill\resizebox{0.47\linewidth}{!}{%
			\begin{tikzpicture}
			\begin{axis}[%
			xmode=normal,
			ymode=log,
			xlabel={$d$},
			xmin=1,xmax=16,
			legend style={at={(0.01,0.01)},anchor=south west},
			]
			\addplot+[no marks,line width=1.0pt] coordinates{
				( 1, 3.02714e-01)
				( 2, 3.08858e-02)
				( 3, 8.64294e-03)
				( 4, 6.77285e-03)
				( 5, 2.67812e-03)
				( 6, 1.44066e-03)
				( 7, 1.05699e-03)
				( 8, 6.01468e-04)
				( 9, 2.33519e-04)
				(10, 1.44330e-04)
				(11, 7.72676e-05)
				(12, 3.70242e-05)
				(13, 2.03524e-05)
				(14, 9.31298e-06)
				(15, 5.01954e-06)
				(16, 2.54304e-06)
				(17, 1.79800e-06)
				(18, 1.49949e-06)
				(19, 8.65168e-07)
				(20, 5.77397e-07)
			}; \addlegendentry{Hankel Singular values};
			\addplot+[] coordinates{
				(4 ,  0.159961 - 0.140850)
				(6 ,  0.148092 - 0.140850)
				(8 ,  0.141219 - 0.140850)
				(10,  0.141069 - 0.140850)
				(12,  0.140850 - 0.140758)
				(14,  0.140850 - 0.140836)
			}; \addlegendentry{$|\mathcal{J}_{d} - \mathcal{J}_{15}|$};
			\end{axis}
			\end{tikzpicture}
		}
	\end{figure}

	In Fig. \ref{fig:fp-Dd} (left) we study the total cost in the LQR stabilised system.
	Here we initialise the Fokker-Planck system with the density function of the uniform distribution on $[-6,6]$.
	We can deduce that an absolute error of about $10^{-3}$ (a relative error of 1\%) is achieved for 1023 points in the initial finite difference discretization.
	Setting $D=1023$ and varying the basis size in the balanced truncation,
	we compare the Hankel singular values and the differences in the total cost in the HJB stabilised system in Fig. \ref{fig:fp-Dd} (right).
	We observe that  $d=10$ dimensions in the reduced model are sufficient to drop the absolute error below the same level of $10^{-3}$.

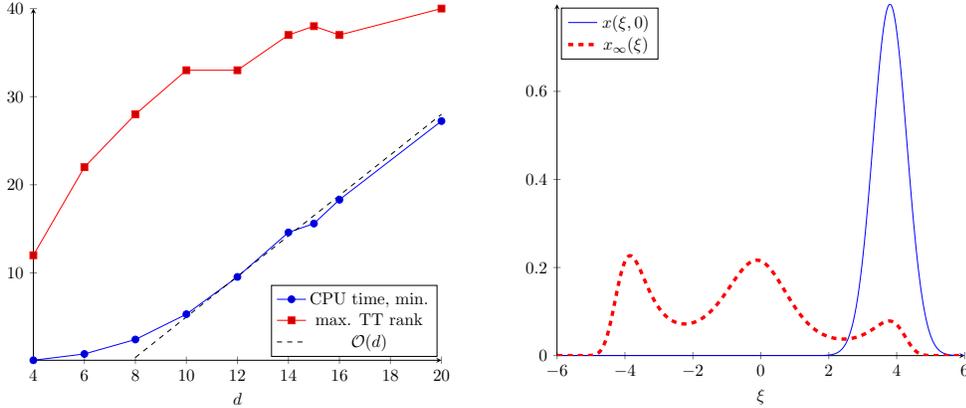
\begin{figure}[h!]
		\centering
		\caption{Left: CPU times and TT ranks for different dimensions~$d$ for the Fokker-Planck problem with the right-sided initial state $x(\xi,0) = \frac{1}{Z}\exp(-2(\xi-3.8)^2)$. Right: initial and equilibrium states.}
		\label{fig:fp-ttimes}
		\resizebox{0.47\linewidth}{!}{%
			\begin{tikzpicture}
			\begin{axis}[%
			xmode=normal,
			ymode=normal,
			xlabel={$d$},
			legend style={at={(0.99,0.01)},anchor=south east},
			]
			\addplot+[] coordinates{
				(4 , 4.986361  /60)
				(6 , 47.882615 /60)
				(8 , 146.279308/60)
				(10, 318.694551/60)
				(12, 572.491130/60)
				(14, 875.573600/60)
				(15, 935.791793/60)
				(16, 1098.19868/60)
				(20, 1633.860  /60)
			}; \addlegendentry{CPU time, min.};
			\addplot+[] coordinates{
				(4 , 12)
				(6 , 22)
				(8 , 28)
				(10, 33)
				(12, 33)
				(14, 37)
				(15, 38)
				(16, 37)
				(20, 40)
			}; \addlegendentry{max. TT rank};
			\addplot+[no marks,dashed,black,domain=8:20] {2.3*x-18}; \addlegendentry{$\mathcal{O}(d)$};
			
			\end{axis}
			\end{tikzpicture}
		}
		\hfill\resizebox{0.47\linewidth}{!}{%
			\begin{tikzpicture}
			\begin{axis}[%
			xmode=normal,
			ymode=normal,
			xlabel={$\xi$},
			xmin=-6,xmax=6,
			legend style={at={(0.01,0.99)},anchor=north west},
			]
			\addplot+[no marks] table[header=true,x=xi,y=x0] {fp-densities.dat}; \addlegendentry{$x(\xi,0)$};
			\addplot+[no marks,line width=2pt,dashed] table[header=true,x=xi,y=xinf] {fp-densities.dat}; \addlegendentry{$x_\infty(\xi)$};
			\end{axis}
			\end{tikzpicture}
		}
	\end{figure}

	In Fig.~\ref{fig:fp-ttimes} (left) we vary the dimension $d$ of the reduced state and investigate CPU times and TT ranks of the value function.
	The TT approximation threshold $\delta=10^{-4}$, the initial shift in Alg. \ref{alg:als} $\mu=5$, and the polynomial degree $n-1=4$.
	We see that the TT ranks stabilize as the dimension increases, and hence the CPU time grows linearly.
	
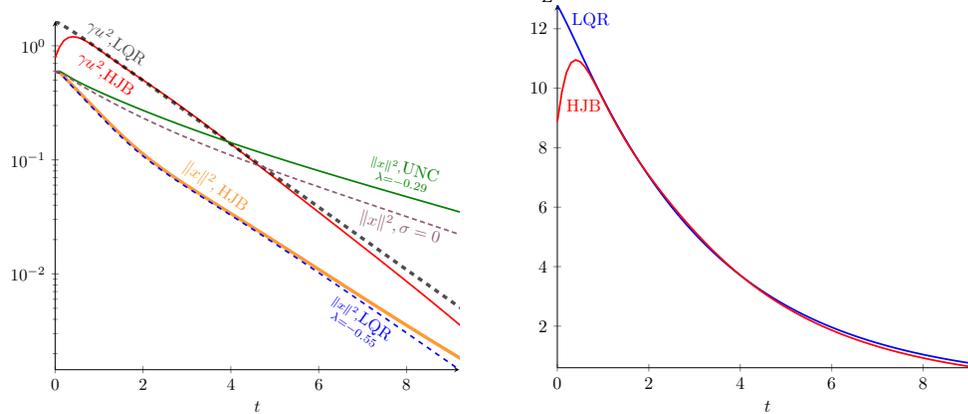
\begin{figure}[h!]
		\centering
		\caption{Running costs (left) and control signals (right) for the reduced Fokker-Planck problem with $d=10$ with the right-sided initial state $x(\xi,0) = \frac{1}{Z}\exp(-2(\xi-3.8)^2)$.}
		\label{fig:fp-cost}
		\resizebox{0.47\linewidth}{!}{%
			\begin{tikzpicture}
			\begin{axis}[%
			xmode=normal,
			ymode=log,
			xmin=0,xmax=9.2,
			xlabel={$t$},
			legend style={at={(0.99,0.99)},anchor=north east},
			]
			
			\addplot+[no marks,line width=2.0pt,orange,opacity=0.8] table[header=true,x=t,y=chjbx]{fp_cost_d10.dat} node[pos=0.15,anchor=south west,rotate=-35]{$\|x\|^2,\mbox{HJB}$};
			\addplot+[no marks,line width=1.0pt] table[header=true,x=t,y=chjbu]{fp_cost_d10.dat} node[pos=0.02,anchor=north west,rotate=-35]{$\gamma u^2$,HJB};
			\addplot+[no marks,line width=1.0pt,dashed,blue] table[header=true,x=t,y=clqrx]{fp_cost_d10.dat} node[pos=0.4,anchor=north east,rotate=-35]{$\substack{\|x\|^2,\mbox{LQR}\\\lambda=-0.55}$};
			\addplot+[no marks,line width=1.0pt,dashed,line width=2pt,opacity=0.7] table[header=true,x=t,y=clqru]{fp_cost_d10.dat} node[pos=0.02,anchor=south west,rotate=-35]{$\gamma u^2$,LQR};
			\addplot+[no marks,line width=1.0pt,green!50!black] table[header=true,x=t,y=cunc]{fp_cost_d10.dat} node[pos=0.35,anchor=south west,rotate=-15]{$\substack{\|x\|^2,\mbox{UNC}\\\lambda=-0.29}$};
			
			\addplot+[no marks,line width=1.0pt,magenta!50!black] table[header=true,x=t,y=chjbx]{fp_cost_d10_ds0.dat} node[pos=0.35,anchor=north west,rotate=-15]{$\|x\|^2,\sigma=0$};
			
			\end{axis}
			\end{tikzpicture}
		}
		\hfill\resizebox{0.47\linewidth}{!}{%
			\begin{tikzpicture}
			\begin{axis}[%
			xmode=normal,
			ymode=normal,
			xmin=0,xmax=9.2,
			xlabel={$t$},
			legend style={at={(0.99,0.99)},anchor=north east},
			]
			
			\addplot+[no marks,line width=1.0pt]
			table[header=true,x=t,y=ulqr]{fp_cost_d10.dat} node[pos=0.02,anchor=west]{LQR};
			\addplot+[no marks,line width=1.0pt] table[header=true,x=t,y=uhjb]{fp_cost_d10.dat} node[pos=0.02,anchor=west]{$\mbox{HJB}$};
			\end{axis}
			\end{tikzpicture}
		}
	\end{figure}

	Moreover, we change the initial distribution to the right-sided state $x(\xi,0) = \frac{1}{Z}\exp(-2(\xi-3.8)^2)$, where $Z = \int_{-6}^{6} \exp(-2(\xi-3.8)^2) d\xi$ is the normalisation constant (Fig~\ref{fig:fp-ttimes}, right).
	It was observed \cite{bkp-hjb-fp-2018} that the free system exhibits a very slow convergence to equilibrium when started from a right-sided distribution, since the particles must flow through a region of low probability.
	In Fig.~\ref{fig:fp-cost} we show the components of the running cost for  the \emph{original} unshifted system,
	both UNControlled and controlled with HJB and LQR signals, obtained for the shifted system.
	We see that the free system converges at a slow rate $\|x\|^2 \sim \exp(-0.29 t)$, while the controller computed for the de-stabilised system can accelerate this rate by almost a factor of 2.
	Note that when the HJB controller is computed for the original system ($\sigma=0$), it accelerates the convergence only a little, so the shift is important to achieve the speedup.
	However, larger shifts make the HJB equation more difficult to solve.
	In particular, for larger shifts $\sigma$ and larger state domain sizes $a$ the stiffness matrix in \eqref{eq:discnd} might become indefinite, and the policy iteration fails to converge.
	The domain size should be large enough to fit the trajectory, e.g.
	for the right-sided initial state the domain size of $a=20$ is necessary to avoid excessive extrapolation of Legendre polynomials.
	This poses certain limitations on the range of possible applications of the TT-HJB approach.
	Nevertheless, when the policy iteration converges the HJB regulator can deliver a lower cost than LQR.

	\subsection*{Conclusion}
	We presented a numerical method for the solution of high-dimensional HJB PDEs arising in optimal feedback control
    for nonlinear dynamical systems. Our algorithm combines a continuous policy iteration  with a tensor-train ansatz for the value function. An important matter of investigation is the identification of a class of optimal control problems where the value function can be accurately represented with a low-rank tensor train structure. For the class of  optimal control problems we have explored in this work, consisting of systems governed by nonlinear parabolic PDEs, we consistently showed that the maximum TT rank in the value function approximation scales linearly with the dimension. This allows us to circumvent the curse of dimensionality up to a great extent, solving HJB PDEs with more than 100 dimensions. Control constraints are effectively enforced through penalties, despite larger TT ranks of the value function.
    The applications of the proposed methodology are extensive. Here we  explored the synthesis of feedback control laws for high-dimensional dynamics arising from the semi-discretisation of nonlinear PDEs. However, high-dimensional dynamics also play a crucial role in aerospace engineering \cite{IFAC}, networks and agent-based models \cite{ACFK}. Finally, our methodology based on spectral approximation  and tensor calculus, opens possibilities for a new error analysis.

	\subsection*{ACKNOWLEDGMENTS} S. Dolgov acknowledges the support of the Engineering and Physical Sciences Research Council through Fellowship EP/M019004/1. K. Kunisch was partially supported by the ERC advanced grant 668998 (OCLOC) under the EU's H2020 research program. Dante Kalise was supported by a public grant as part of the Investissement d’avenir project,reference ANR-11-LABX-0056-LMH, LabEx LMH.

	\bibliographystyle{siam}
	\bibliography{references}
\end{document}